\documentclass[12pt]{amsart}
\usepackage{amsfonts}
\usepackage{amsmath}
\usepackage{amsxtra}
\usepackage{amssymb,latexsym}
\usepackage[mathcal]{eucal}
\usepackage{amscd}

\makeindex

\newtheorem{theo}{{\bfseries Theorem}}[section]
\newtheorem{prop}[theo]{{\bfseries Proposition}}
\newtheorem{lem}[theo]{{\bfseries Lemma}}
\newtheorem{cor}[theo]{{\bfseries Corollary}}
\newtheorem{df}[theo]{{\bfseries Definition}}

\newtheorem{ex}[theo]{{\bfseries Example}}

\newtheorem{ques}[theo]{{\bfseries Question}}

\def \ol {\overline}
\def \N {\mathbb N}

\def \Z {\mathbb Z}
\def \R {\mathbb R}

\def \A {\mathcal A}
\def \B {\mathcal B}
\def \CC {\mathcal C}

\def \NN {\mathcal N}

\def \RR {\mathcal R}

\def \G {\mathcal G}

\def \ep {\epsilon}

\def \d {\delta}

\usepackage{amssymb,latexsym}
\usepackage[mathcal]{eucal}
\usepackage{amscd}
\usepackage{amsmath}
\usepackage{graphicx}
\usepackage{amsfonts}
\usepackage{amscd}

\numberwithin{equation}{section}

\begin{document}

\title{\bfseries  Varieties of Mixing}
\vspace{1cm}
\author{Ethan Akin and Jim Wiseman}
\address{Mathematics Department \\
    The City College \\ 137 Street and Convent Avenue \\
       New York City, NY 10031, USA     }
\email{ethanakin@earthlink.net}

\address{Department of Mathematics \\
Agnes Scott College \\ 141 East College Avenue \\ Decatur, GA 30030, USA }
\email{jwiseman@agnesscott.edu}

\thanks{The second author was supported by a grant from the Simons Foundation (282398, JW)}

\date{October, 2017}

\begin{abstract} We consider extensions of the notion of topological transitivity for a dynamical system $(X,f)$.  
In addition to chain transitivity, we define
strong chain transitivity and vague transitivity. Associated with each there is a notion of mixing, defined by transitivity of
the product system $(X \times X, f \times f)$. These extend the concept of weak mixing which is associated with topological transitivity.
Using the barrier functions of Fathi and Pageault, we obtain for each of these extended notions a dichotomy result that a transitive
system of each type either satisfies the corresponding mixing condition or else factors onto an appropriate type of equicontinuous minimal
system. The classical dichotomy result for minimal systems follows when it is shown that a minimal system is weak mixing if and only if it is
vague mixing.
\end{abstract}


\vspace{.5cm} \maketitle

\tableofcontents

\section{Introduction}

Consider a \emph{dynamical system} consisting of a continuous map $f$ on a
compact metrizable space $X$. We call the system
\emph{topologically transitive } when for some $x \in X$
the orbit $\{ f(x), f^{2}(x), \dots \}$ is dense in $X$. Such a \emph{transitive point} is necessarily recurrent.
For a topologically transitive system the set of such transitive points is an $f$ invariant,
dense $G_{\d}$ subset and the map is surjective.
The system is \emph{minimal} when every
orbit is dense.  It is \emph{weak mixing} when $f \times f$ is topologically transitive on $X \times X$.

When $f$ is a minimal homeomorphism, either $f$ is weak mixing or it has a nontrivial equicontinuous factor.

Our purpose here is to consider other forms of transitivity and the related notions of
mixing and to obtain for them analogous dichotomy results.

If $[x_1,\dots,x_n]$ is a sequence in a compact metric space $(X,d)$ and $x,y \in X$ then for the dynamical system $f$ on $X$ we define
the $xy$ \emph{chain-bound} to be $\max[d(x,x_1), d(f(x_1),x_2),\dots, d(f(x_{n-1}),x_n), d(f(x_n),y)]$\\ and the  \emph{chain-length}
to be $d(x,x_1) + d(f(x_1),x_2) + \dots +  d(f(x_{n-1}),x_n) + d(f(x_n),y)$. Thus, we begin with $x_1$ near $x$, iterate $n$ times, terminating
at $f(x_n)$ near $y$.  At each step we make an error measured by $d(f(x_i),x_{i+1})$. The chain-bound and chain-length are alternative
ways of measuring the total error.

We define the \emph{chain relation} $\CC_{d} f $ to be the set of pairs $(x,y) \in X \times X$ such that for every $\ep > 0$ there exists
$[x_1,\dots,x_n]$ with $xy$ chain-bound less than $\ep$. The \emph{strong chain relation} $\A_{d} f $ is
the set of pairs $(x,y) \in X \times X$ such that for every $\ep > 0$ there exists
$[x_1,\dots,x_n]$ with $xy$ chain-length less than $\ep$. Each of these is a closed, transitive relation which contains $f$. From uniform continuity
it is clear that the chain relation $\CC_d f$ is independent of the choice of admissible metric $d$ on the compact metrizable space $X$. The strong
chain relation $\A_d f$ does depend upon the metric. On the other hand, by intersecting  the relations $\A_d f$ with $d$ varying over
all admissible metrics, we obtain $\G f$, the smallest closed, transitive relation which contains $f$.

We call $f$ \emph{chain transitive} when $\CC_d f = X \times X$, \emph{strong chain transitive} when $\A_d f = X \times X$ and \emph{vague transitive} when
$\G f = X \times X$.  The map is called \emph{chain mixing, strong chain mixing} or \emph{vague mixing} when the product function $f \times f$ on the product
metric space $(X \times X, d \times d)$ is chain transitive, strong chain transitive or vague transitive.

\begin{theo}\label{theoi00} Let $f$ be a continuous map on a compact metric space $(X,d)$.
\begin{enumerate}
\item[(a)] If $f$ is chain transitive, then either $f$ is chain mixing or there exists a continuous
function mapping $f$ onto a non-trivial periodic orbit.
\item[(b)] If $f$ is strong chain transitive, then either $f$ is strong chain mixing or there exists a Lipschitz
function mapping $f$ onto a non-trivial minimal, isometric homeomorphism.
\item[(c)] If $f$ is vague transitive, then either $f$ is vague mixing or there exists a continuous
function mapping $f$ onto a non-trivial minimal, equicontinuous homeomorphism.
\end{enumerate}
In each case the ``or'' is exclusive.
\end{theo}

The \emph{non-wandering relation} $\NN f$ consists of those pairs $(x,y) \in X \times X$ such that for every $\ep > 0$ there exists
$x_1 \in X$ and $n \in \N$ with $d(x,x_1), d(f^n(x_1),y) < \ep$. So for this relation, errors occur only at the beginning and end of the
orbit sequence. The system $f$ is topologically transitive exactly when $\NN f = X \times X$. However, because $\NN f$ is not in general
a transitive relation, we do not obtain the sort of dichotomy result as in Theorem \ref{theoi00}. If $f$ is a minimal homeomorphism then we do, because
of the following.

\begin{theo}\label{theoi00a} If $f$ is a minimal homeomorphism on a compact metric space $(X,d)$, then $f$ is weak mixing iff it is
vague mixing. \end{theo}

It will be convenient to use the language of relations following \cite{A93}. All our spaces
are compact, metrizable spaces.

A relation  $f : X \to Y$ is  a subset of $X \times Y$ with
$f(x) = \{ y \in Y : (x,y) \in f \}$ for $x \in X$,
and let $f(A) = \bigcup_{x \in A} \ f(x)$ for $A \subset X$.
So $f$ is a mapping when $f(x)$ is a singleton set for every $x \in X$, in which case we will use the notation
$f(x)$ for both the singleton set and the point contained therein. For
example, the identity map on  $X$ is $1_X = \{ (x,x) : x \in X \}$.
We call $f$  a closed relation when
it is a closed subset of $X \times Y$ with the product
topology.

All our relations will be assumed to be nonempty.

For a relation $f : X \to Y$ the \emph{inverse relation}
$f^{-1} : Y \to X $ is $ \{ (y,x) : (x,y) \in f \}$.
Thus, for $B \subset Y$, $f^{-1}(B) = \{ x : f(x) \cap B \not= \emptyset \}.$

If $f : X \to Y$ and $g : Y \to Z$ are relations then
the \emph{composition} $g \circ f : X \to Z $ is
$\{ (x,z) : $ there exists $y \in Y$ such that $(x,y) \in f$ and $(y,z) \in  g \}$.  That is, $g \circ f$ is the image
of $(f \times Z) \cap (X \times g)$ under the projection $\pi_{13} : X \times Y \times Z \to X \times Z$.
As with maps, composition of relations is clearly associative. The composition of closed relations is closed.

 The \emph{domain} of a relation $f : X \to Y$ is
\begin{equation}\label{introeq01}
Dom(f) \ = \ \{ x : f(x) \not= \emptyset \} \ = \ f^{-1}(Y).
\end{equation}
We call a relation \emph{surjective} if $Dom(f) = X$ and $Dom(f^{-1}) = Y$, i.e.\ $f(X) = Y$ and $f^{-1}(Y) = X$.

If $f_1 : X_1 \to Y_1$ and $f_2 : X_2 \to Y_2$ are relations, then the \emph{product relation}
$f_1 \times f_2 : X_1 \times X_2 \to Y_1 \times Y_2$ is
$\{ ((x_1,x_2),(y_1,y_2)) : (x_1,y_1) \in f_1, (x_2,y_2) \in f_2 \}$.

We call $f$ a \emph{relation on $X$} when $X = Y$.  In that case, we define, for $n \geq 1$
$f^{n+1} = f \circ f^n = f^n \circ f$ with $f^1 = f$.  By definition, $f^0 = 1_X$ and $f^{-n} = (f^{-1})^n$.
If $A \subset X$, then
$A$ is called $f$ \emph{$^+$invariant} if $f(A) \subset A$ and $f$ \emph{invariant} if $f(A) = A$.

For a pseudo-metric $d$ on $X$, we define the relations on $X$
\begin{equation}\label{introeq02}
V^d_{\ep} = \{ (x,y) : d(x,y) < \ep \}, \quad \bar V^d_{\ep} = \{ (x,y) : d(x,y) \leq \ep \}.
\end{equation}

A relation $f$ on $X$ is \emph{reflexive} when $1_X \subset f$, \emph{symmetric} when
$f^{-1} = f$ and \emph{ transitive} when $f \circ f \subset f$.

Given a closed relation $f$ on $X$, the \emph{orbit closure relation} $\RR f$ is defined by
$\RR f(x) = \ol{ \bigcup_{n \in \N} f^n(x) }$. This is usually a proper subset of
$\NN f = \ol{ \bigcup_{n \in \N} f^n }$. The latter is a closed relation but not usually transitive.
We define $\G f$ to be the smallest closed, transitive relation which contains $f$.

For a closed relation $f$ on a compact metric space $(X,d)$ we will also define the chain
relation $\CC_d f$ and the strong chain relation $\A_d f$.
These are closed transitive relations with $\NN f \subset \G f \subset \A_d f \subset \CC_d f$.

We will call a relation $f$ vague transitive
when $\G f = X \times X$, strong chain transitive when $\A_d f = X \times X$ and
chain transitive when $\CC_d f = X \times X$.
For each of these there is a corresponding notion of mixing which is transitivity of $f \times f$ on $X \times X$.

By using the barrier functions of Easton \cite{E}, Pageault \cite{P} and Fathi \cite{FP},
we are able to get dichotomy results for these
notions of transitivity and mixing, see Theorem \ref{theoi00}.

We will need some simple results about pseudo-metrics.

\begin{prop}\label{propi01} Let $X$ be a compact space.
\begin{enumerate}
\item[(a)] Let $\rho$ be a pseudo-metric on $X$. The following are equivalent
and when they hold we call $\rho$ a continuous pseudo-metric
on $X$.
\begin{itemize}
\item[(i)] The identity map from $X$ to the pseudo-metric space $(X,\rho)$ is continuous.
\item[(ii)] The map $\rho : X \times X \to \R$ is continuous.
\item[(iii)] For every $x \in X$ the map $y \mapsto \rho(x,y)$ is continuous.
\item[(iv)] For every $\ep > 0$ the set $V^{\rho}_{\ep}$ is an open subset of $X \times X$.
\item[(v)]  For every $\ep > 0$ the set $\bar V^{\rho}_{\ep}$ is a neighborhood of the diagonal $1_X \subset X \times X$.
\item[(vi)] For every $\ep > 0$  and $x \in X$  the ball $\bar V^{\rho}_{\ep}(x)$ is a neighborhood of $x$ in $X$.
\end{itemize}

\item[(b)] If $d$ is a continuous metric on $X$ then the topology of $X$ is that of the metric space $(X,d)$, i.e. $d$ is an
\emph{admissible metric} for the space $X$.

\item[(c)] If $X$ is metrizable then a pseudo-metric $\rho$ on $X$ is continuous
iff there exists an admissible metric $d$ on $X$ such that
$d \geq \rho$.
\end{enumerate}
\end{prop}

{\bfseries Proof:} (a): (i) $\Rightarrow$ (ii): On the product pseudo-metric space
$(X,\rho) \times (X,\rho)$ the map $\rho$ is continuous.
Compose with the continuous map from $X \times X$ to $(X,\rho) \times (X,\rho)$.

(ii) $\Rightarrow$ (iii) $\Rightarrow$ (vi) and (ii) $\Rightarrow$ (iv)
$\Rightarrow$ (v) $\Rightarrow$ (vi) are obvious.

Finally, it is clear that (vi) $\Rightarrow$ (i).

(b): If $d$ is a continuous metric then $X \to (X,d)$ is a continuous bijection
from a compact space to a Hausdorff space and so is a
homeomorphism.

(c): If $d \geq \rho$ then  $\bar V^{\rho}_{\ep} \supset \bar V^{d}_{\ep}$ and so
$\rho$ satisfies (v) of (a). Conversely, if $\rho$ is
a continuous pseudo-metric and $d_1$ is a continuous metric then $d = d_1 + \rho$ is a continuous metric.

$\Box$ \vspace{.5cm}

A pseudo-metric (or metric) $d$ is a \emph{pseudo-ultrametric} (or an \emph{ultrametric})
when it satisfies the ultrametric strengthening of
the triangle inequality: $d(x,z) \leq \max(d(x,y),d(y,z))$. This is equivalent
to the assumption that  $V^d_{\ep}$ (or that
 $\bar V_{\ep}^d$) is an equivalence relation for every $\ep > 0$. Note that if
 $E$ is an equivalence relation on a compact space $X$
 and $E$ is a neighborhood of the diagonal $1_X$ then each equivalence class $E(x)$
 is a neighborhood of each of its points and so is open. The
 complement of $E(x)$ is a union of equivalence classes and so $E(x)$ is clopen.
 By compactness there are only finitely many equivalence
 classes and so $E = \bigcup_x E(x) \times E(x)$ is clopen in $X \times X$.

 It follows that if a compact space admits a continuous ultrametric, then the clopen sets
 form a basis for the topology. Such a space
 is called \emph{zero-dimensional}. Conversely, if $\B$ is the set of clopen subsets
 of a compact metrizable space $X$ then
 $\B$ is a countable set (as each $B \in \B$ is a finite union of some members of a countable basis).
 If $\B$ forms a basis
 then $j : X \to \{ 0, 1 \}^{\B}$ defined by $j(x)_B = 1 \ \Leftrightarrow \ x \in B$ is an embedding.
 If $d$ is the zero-one metric on $\{ 0,1 \}$
 then $d(a,b) = \max_{n}  2^{-n}d(a_{n},b_{n})$ is an ultrametric on $\{ 0, 1 \}^{\N}$. Hence,
 a compact, zero-dimensional metrizable space
 admits a continuous ultrametric.

 \vspace{1cm}

\section{Chains}
\vspace{.5cm}

Let $(X,d)$ be a compact metric space and $f$ be a closed  relation on $X$.
An $f$-chain  of length $n$ is an element $C $ of
the $n$-fold product $f^{\times n}$ with $n \geq 1$. So $C
 = [(a_1,b_1), \dots,(a_n, b_n)]$ with $n \geq 1$,
such that $(a_i,b_i) \in f$ for $i = 1, \dots, n$. .
In particular, if $f$ is a map then $b_i = f(a_i)$ for $i = 1, \dots, n$.
 The \emph{length} of $C$ is $\# C = n$, the
\emph{$xy$ chain-bound}  of $C$ is $|xCy| = max_{i=1}^{n+1} \{ d(b_{i-1},a_i) \}$,
where we let $b_0 = x$ and $a_{n+1} = y$. The \emph{$xy$ chain-length}
 of $C$ is
$||xCy|| = \sum_{i=1}^{n+1} \{ d(b_{i-1},a_i) \}$. Clearly, $|xCy| \leq ||xCy||$.

Clearly, with $x,y,z \in X$
\begin{align}\label{c01}
\begin{split}
 |xCz| \ \leq \ |xCy| + d(y,z), \quad &|zCy| \ \leq \ d(z,x) + |xCy|, \\
  ||xCz|| \ \leq \ ||xCy|| + d(y,z), \quad &||zCy|| \ \leq \ d(z,x) + ||xCy||
\end{split}
\end{align}

If $C \in f^{\times n}$ and $D \in f^{\times m}$, then the
concatenation $C \cdot D \in f^{\times (n + m)}$ so that $\# C \cdot D = \# C + \# D$.
With $x, y, z \in X$
\begin{align}\label{c02}
\begin{split}
 |xC \cdot Dy| \ \leq \ &|xCz| + |zDy| \\ ||xC \cdot Dy|| \ \leq \ &||xCz|| + ||zDy||.
\end{split}
\end{align}
In the case of the chain-bound, the only sum estimates occur between $C$ and $D$. Hence, if
$C \in f^{n}, D \in f^{m}, E \in f^{p}$ then for the concatenation $C \cdot D \cdot E$ we have
for $x,y,z,w \in X$:
\begin{equation}\label{c02a}
|x C \cdot D \cdot E y| \ \leq \ \max( |x C z| + |z D w|, |z D w| + |w E y| ).
\end{equation}

If $C = [(a_1,b_1),\dots,(a_n,b_n)] \in f^{\times n}$ we let
$C^{-1} = [(b_n,a_n), \dots, (b_1,a_1)] \in (f^{-1})^{\times n}$.
Clearly,
\begin{equation}\label{c02b}
||x C y|| = ||y C^{-1} x||, \quad \text{and} \quad |x C y| = |y C^{-1} x|.
\end{equation}

If $h : (X,d) \to (X_1, d_1)$ is a continuous map of compact metric spaces we define, for $\d > 0$,
$o(h,\d) = \sup \{ d_1(h(x),h(y)) : (x,y) \in \bar V^d_{\d} \}$. Uniform
continuity says that for every $\ep > 0$ there exists
$\d > 0$ such that $o(h,\d) < \ep$. If $h$ is Lipschitz with Lipschitz constant $L$ then $o(h,\d) \leq L \d$.

We say that $h$ maps $f$ on $X$ to $f_1$ on $X_1$ if $(h \times h)(f) \subset f_1$.
Note that $(h \times h)(f) = h \circ f \circ h^{-1}$.

If $h$ maps $f$ to $f_1$ and $C  \in f^{\times n }$  then applying $h$ to
each term of the sequence we obtain $h(C) \in f_1^{\times n}$
with $\# h(C) = \# C$.  For $x,y \in X$
\begin{equation}\label{c03}
 |h(x)h(C)h(y)| \ \leq \ o(h,|xCy|),
\end{equation}
and if $h$ has Lipschitz constant $L$ then
\begin{equation}\label{c04}
|h(x)h(C)h(y)| \ \leq \ L|xCy|, \qquad ||h(x)h(C)h(y)|| \ \leq \ L||xCy||.
\end{equation}

We define the \emph{barrier functions} $m^f_d, \ell^f_d : X \times X \to \R$ by
\begin{align}\label{c05}
\begin{split}
m^f_d(x,y) \ = \ &inf \{ |xCy| : C \in f^{\times n}, n \in \N \}, \\
\ell^f_d(x,y) \ = \ &inf \{ ||xCy|| : C \in f^{\times n}, n \in \N \}.
\end{split}
\end{align}
Obviously, $m^f_d \leq \ell^f_d$. If $(x,y) \in f$ then $C = [(x,y)] \in f^{\times 1}$ with $||xCy || = 0$. Hence,
\begin{equation}\label{c05a}
(x,y) \in f \quad \Longrightarrow \quad m^f_d(x,y) = \ell^f_d(x,y) = 0.
\end{equation}

From (\ref{c01}) we see that for $x,y,z \in X$
\begin{align}\label{c06}
\begin{split}
|m^f_d(x,y) - m^f_d(x,z)| \ \leq \ d(y,z), \quad &|m^f_d(x,y) - m^f_d(z,y)| \ \leq \ d(z,x), \\
|\ell^f_d(x,y) - \ell^f_d(x,z)| \ \leq \ d(y,z), \quad &|\ell^f_d(x,y) - \ell^f_d(z,y)| \ \leq \ d(z,x).
\end{split}
\end{align}

In particular, using $z = y$, we obtain
\begin{align}\label{c06a}
\begin{split}
m_f^d(y,y) = 0 \qquad &\Longrightarrow \qquad m_f^d(x,y) \ \leq \ d(x,y),\\
\ell_f^d(y,y) = 0 \qquad &\Longrightarrow \qquad \ell_f^d(x,y) \ \leq \ d(x,y).
\end{split}
\end{align}

From (\ref{c02}) we see that
\begin{align}\label{c07}
\begin{split}
m^f_d(x,z) \ \leq \ \ &m^f_d(x,y) \ + \ m^f_d(y,z), \\
 \ell^f_d(x,z) \ \leq \ \ &\ell^f_d(x,y) \ + \ \ell^f_d(y,z).
\end{split}
\end{align}
Furthermore, from (\ref{c02a}) we obtain
\begin{equation}\label{c07a}
m^f_d(x,z) \ \leq \ \max[m^f_d(x,y) + m^f_d(y,y), m^f_d(y,y) + m^f_d(y,z)].
\end{equation}

From (\ref{c02b}) we see that
\begin{equation}\label{c07b}
m^f_d(x,y) \ = \ m^{f^{-1}}_d(y,x) \quad \text{and} \quad \ell^f_d(x,y) \ = \ \ell^{f^{-1}}_d(y,x).
\end{equation}

If $h : (X,d) \to (X_1, d_1)$ is  continuous, mapping $f$ to $f_1$ then
\begin{equation}\label{c08}
m^{f_1}_{d_1}(h(x),h(y)) \ \leq \ o(h,m^f_d(x,y)), \hspace{3cm}
\end{equation}
and if $h$ has Lipschitz constant $L$ then
\begin{align}\label{c09}
\begin{split}
m^{f_1}_{d_1}(h(x),h(y)) \ &\leq \ L \cdot m^f_d(x,y), \\ \ell^{f_1}_{d_1}(h(x),h(y)) \ &\leq \ L \cdot \ell^f_d(x,y).
\end{split}
\end{align}

For a closed relation $f$ on a compact metric space $(X,d)$ we define the
chain relation $\CC_d f$ and the strong chain relation $\A_d f$
to be the relations on $X$ given by
\begin{equation}\label{c10}
\CC_d f = \{ (x,y) : m^f_d(x,y) = 0 \}, \qquad \A_d f = \{ (x,y) : \ell^f_d(x,y) = 0 \}.
\end{equation}

From $m^f_d \leq \ell^f_d$ and (\ref{c05a}) we obtain
\begin{equation}\label{c11}
f \ \subset \ \A_d f \ \subset \ \CC_d f. \hspace{2cm}
\end{equation}

The functions $\ell^f_d$ and $m^f_d$ are continuous and so $\A_d f$ and $\CC_d f$
are closed relations. From the triangle inequalities
(\ref{c07}) it follows that the relations are transitive. From (\ref{c07b}) we
see that  $\A_d (f^{-1}) = (\A_d f)^{-1}$ and
$\CC_d (f^{-1}) = (\CC_d f)^{-1}$.  So we may omit the parentheses.

From (\ref{c06a}) we obtain
\begin{align}\label{c11a}
\begin{split}
1_X \ \subset \ \CC_d f \qquad &\Longrightarrow \qquad m_f^d(x,y) \ \leq \ d(x,y),\\
1_X \ \subset \ \A_d f \qquad &\Longrightarrow \qquad \ell_f^d(x,y) \ \leq \ d(x,y),
\end{split}
\end{align}
for all $x,y \in X$.

If $h : (X,d) \to (X_1, d_1)$ is  continuous, mapping $f$ to $f_1$, then (\ref{c08})
implies that $h$ maps $\CC_d f$ to $\CC_{d_1} f_1$
and if $h$ is Lipschitz then by (\ref{c09}) it maps $\A_d f$ to $\A_{d_1} f_1$.
In particular, it follows that $\CC_d f$ is independent of
the choice of metric on $X$. On the other hand, the strong chain relation $\A_d f$
does depend upon the metric.  Notice that if
$d_1, d_2$ are metrics on $X$ then $d = \max(d_1,d_2)$ is a metric on $X$ and
the identity maps from $(X,d)$ to $(X,d_1)$ and $(X,d_2)$ have
Lipschitz constant $1$. It follows that
\begin{equation}\label{c12}
\A_d f \ \subset \ (\A_{d_1} f) \cap (\A_{d_2} f). \hspace{2cm}
\end{equation}
So as $d$ varies over the set of admissible metrics on $X$, the
collection of relations $\{ \A_d f \}$ is a filter-base. Theorem 5.14 and Proposition
6.15 of \cite{AW17}, which extend a theorem of \cite{FP}, imply that

\begin{theo}\label{theoc01} If $f$ is a closed relation on a compact
metrizable space $X$ then $\G f = \bigcap_d  \ \{ \A_d f \}$ with
$d$ varying over the set of admissible metrics on $X$. \end{theo}

$\Box$ \vspace{.5cm}

The original barrier functions of  \cite{P} and \cite{FP} follow \cite{C} and \cite{E} in
not allowing an initial jump and defining $M^f_d, L^f_d : X \times X \to \R$ by
\begin{align}\label{c05d}
\begin{split}
M^f_d(x,y)  =  &\inf  \{ |xCy| : C = [(a_1,b_1),\dots,] \in f^{\times n}, n \in \N, \  a_1 = x \}, \\
L^f_d(x,y)  =  &\inf  \{ ||xCy|| : C = [(a_1,b_1),\dots,] \in f^{\times n}, n \in \N, \  a_1 = x \}.
\end{split}
\end{align}

When $f$ is a continuous map the alternative definitions yield similar results.

\begin{lem}\label{lemc01a} Let $f$ be a continuous map on a compact
metric space $(X,d)$ and $x,y \in X$. For every $\ep > 0$, let $\ep/2 > \d > 0$ be such that
$o(f,\d) \leq \ep/2$. For $C = [(a_1,f(a_1)),(a_2,f(a_2)),\dots, (a_n,f(a_n)] \in f^{\times n}$
let $\tilde C = [(x,f(x)),(a_2,f(a_2)),\dots, (a_n,f(a_n)] \in f^{\times n}$.
\begin{align}\label{c05e}
\begin{split}
|x C y| < \d \quad &\Longrightarrow \quad |x \tilde C y| < \ep, \\
||x C y|| < \d \quad &\Longrightarrow \quad ||x \tilde C y|| < \ep.
\end{split}
\end{align}
\end{lem}

{\bfseries Proof:} If $|x C y| < \d $ then $d(x,a_1) < \d$ and so
\begin{equation}\label{c05f}
d(f(x),a_2) \leq d(f(x),f(a_1)) + d(f(a_1),a_2) \leq \ep/2 + d(f(a_1),a_2).
\end{equation}
Hence, $|x \tilde C y| \leq |x C y| + \ep/2$ and  $||x \tilde C y|| \leq ||x C y|| + \ep/2$.

$\Box$ \vspace{.5cm}

\begin{cor}\label{corc01b} If $f$ is a continuous map on a compact
metric space $(X,d)$, then
\begin{equation}\label{c10b}
\CC_d f = \{ (x,y) : M^f_d(x,y) = 0 \}, \qquad \A_d f = \{ (x,y) : L^f_d(x,y) = 0 \}.
\end{equation}
\end{cor}

{\bfseries Proof:} Obvious from Lemma \ref{lemc01a}.

$\Box$ \vspace{.5cm}

\begin{df}\label{defc02} A relation $f$ on a pseudo-metric space $(X,\rho)$
is called an \emph{isometry} when it is a  relation
such that
\begin{equation}\label{c13}
(x_1,y_1), (x_2,y_2) \in f  \quad \Longrightarrow   \quad \rho(x_1,x_2) = \rho(y_1,y_2).
\end{equation}
\end{df}
\vspace{.5cm}

In particular, if $f$ is a  map then it is an isometry for $(X,\rho)$ when $\rho(x_1,x_2) = \rho(f(x_1),f(x_2))$.

For a map $f$ on a compact metric space $(X,d)$ a point $x$ is an \emph{equicontinuity point} when for every
neighborhood $\ep > 0$  there exists $\d > 0$
such that $d(x,y) < \d$ implies
$d(f^n(x),f^n(y)) \leq \ep$ for $n = 0,1,\dots$.
Associated with $d$ and $f$ we define the metric $d_f$ by
\begin{equation}\label{c13a}
d_f(x,y) \ = \ \sup \{ d(f^n(x),f^n(y)) : n = 0,1,\dots \}. \hspace{2cm}
\end{equation}
The point $x$ is an equicontinuity point when it has neighborhoods of arbitrarily small $d_f$ diameter.

When every point is an equicontinuity point, then $f$ is called an \emph{equicontinuous map} and
the $\d$ above can be chosen independent of $x$.

\begin{theo}\label{theoc03} Let $(X,d)$ be a compact metric space.
\begin{enumerate}
\item[(a)] A map $f$ on $(X,d)$ is equicontinuous iff $d_f$ is an admissible  metric on $X$.

\item[(b)] If a relation $f$ is an isometry on $(X,d)$ with $Dom(f) = X$,
then $f$ is a homeomorphism and it is equicontinuous with $d_f = d$.

\item[(c)] If $f$ is a surjective, equicontinuous map on $(X,d)$, then
 $f$ is a homeomorphism on $X$ which is an
isometry on $(X,d_f)$.
Furthermore, every point $x \in X$ is recurrent for $f$, i.e. \ $x \in \RR f(x)$ and each
$\RR f(x)$ is a minimal invariant subset.

\item[(d)] If a map $f$ is an isometry on $(X,d)$ then $\A_d f(x) = \RR f(x)$
for all $x \in X$. If, in addition, $d$ is an ultrametric, then
$\CC_d f(x) = \RR f(x)$ for all $x \in X$.
\end{enumerate}
\end{theo}

{\bfseries Proof:} (a): Obvious.

(b): If $(x,y_1), (x,y_2) \in f$ then $d(y_1,y_2) = d(x,x) = 0$ and so $y_1 = y_2$ since $d$ is a metric. That is,
$f$ is a map. Clearly, $d_f = d$ and so $f$ is equicontinuous. Since $f$ is an isometric map, it is clearly injective.
Let $A = \bigcap_{n \in \N} f^n(X)$ so that $A$ is $f$ invariant. If $\ep > 0$ there exists $n \in \N$ so that $f^n(X) \subset V_{\ep}(A)$.
Hence, with $d(x,A) = \inf \{ d(x,y) : y \in A \}$ for $x \in X$, we see that $d(f^n(x),A) \to 0$ as $n \to \infty$. But since
$f$ is an isometry and $A$ is invariant, $d(f^n(x),A)$ does not vary with $n$. Thus, $d(x,A) = 0$ for all $x \in X$ and so $f$ is surjective.
Thus, $f$ is a homeomorphism.

(c): By \cite{A96} Proposition 2.4, $f$ is a homeomorphism which is an
isometry on $(X,d_f)$ and with every point recurrent.
Hence, the restriction of $f$ to $\RR f(x)$ is an equicontinuous, topologically
transitive map and so it is minimal by \cite{AuY80} Theorem 4.

(d): In any case, $\RR f \subset \G f \subset \A_d f \subset \CC_d f$. Assume that $y \in \A_d f(x)$.
Given $\ep > 0$ there exists $C = [(a_1,b_1),\dots,(a_n,b_n)]$
such that $||xCy|| < \ep$. Since $f$ is a map, $b_i = f(a_i)$ for $i = 1, \dots, n$. With $b_0 = x, a_{n+1} = y$ let
$\ep_i = d(b_i,a_{i+1}) = d(f(a_i),a_{i+1})$
for $i = 0, \dots, n$. Thus, $\sum_i \ep_i < \ep$. Because $f$ is an isometry,
\begin{equation}\label{c14}
\begin{split}
\ep_0 = d(f^n(x),f^n(a_1)), \ep_1 = d(f^n(a_1),f^{n-1}(a_2)),  \dots, \\
 \ep_{n-2} = d(f^2(a_{n-1}),f(a_n)), \ep_n = d(f(a_n),y).
\end{split}
\end{equation}
Hence, $d(f^n(x),y) \leq \sum_i \ep_i < \ep$. Since, $\ep > 0$ was arbitrary, it follows that $y \in \RR f(x)$.

Now assume that $d$ is an ultrametric and $\ep > 0$. $\bar V^d_{\ep}$ is a
clopen equivalence relation. Hence, for
any subset $A \subset X$, $\bar V^d_{\ep}(A)$ is a union of a finite number
of clopen equivalence classes and so is clopen.
Since $f$ is an isometry, $\bar V^d_{\ep}$ is
$f \times f$ $^+$invariant. Because $\RR f(x)$ is $f$ $^+$invariant,
the clopen set $K = \bar V^d_{\ep}(\RR f(x))$ is $^+$invariant.
If $\d > 0$ is smaller than the distance from $K$ to $X \setminus K$,
then if $z \in K$ and $C \in f^{\times n}$ satisfies
$|zCy| < \d$ then, inductively, $a_i \in K$ and so $b_i = f(a_i) \in K$.
Hence, $y \in K$. In particular, $\CC_d f(x) \subset K$.
As $\ep > 0$ was arbitrary, $\CC_d f(x) \subset \RR f(x)$.

$\Box$ \vspace{.5cm}

{\bfseries Remark:} An \emph{odometer} is an inverse limit of a sequence of periodic
orbits of increasing length. A map is
an odometer iff it is a minimal, equicontinuous homeomorphism on
an infinite, compact, zero-dimensional metrizable space, see, e.g. \ \cite{A96} Theorem 3.5.
In general a minimal map on a finite set consists of a single periodic orbit.
\vspace{.5cm}

A pseudo-metric space $(X,\rho)$ has a metric space quotient $(\bar X, \bar \rho)$
obtained as the space of equivalence classes
with respect to zero-set of $\rho$,
\begin{equation}\label{c15}
Z_{\rho} \ = \ \{ (x_1,x_2) : \rho(x_1,x_2) = 0 \}. \hspace{3cm}
\end{equation}
Using the quotient map $\pi : (X,\rho) \to (\bar X, \bar \rho)$
the metric $\bar \rho$ is well-defined by $\bar \rho(\pi(x_1),\pi(x_2)) = \rho(x_1,x_2)$.

\begin{prop}\label{propc04} Let $f$ be a surjective, closed relation on a
compact metrizable space $X$ and let $\rho$ be a continuous pseudo-metric on
$X$. The quotient space $(\bar X, \bar \rho)$ is a compact
metric space with $\pi : (X,\rho) \to (\bar X, \bar \rho)$ Lipschitz with Lipschitz constant 1.

 If $f$ is an isometry on $(X,\rho)$ then
$\bar f = (\pi \times \pi)(f)$ is an isometry on $(\bar X, \bar \rho)$.
The relation $\bar f$
 is an equicontinuous homeomorphism
on $\bar X$.\end{prop}

{\bfseries Proof:}  It is clear that $\bar f$ is an isometry on $(\bar X, \bar \rho)$ and
so by \ref{theoc03} it is an equicontinuous
homeomorphism.

$\Box$ \vspace{1cm}

\section{Mixing}
\vspace{.5cm}


We begin by noting that \cite{AW17} Proposition 6.10 and 6.11 imply that $f$ is
necessarily a surjective relation if $\CC_d f = X \times X$, and so, a fortiori, if $\A_d f = X \times X$
or $\G f = X \times X$.

If $f$ is a closed relation on a compact metric space $(X,d)$ then we define the product relation $f^{(n)}$ on the
metric space $(X^{(n)} , d^{(n)})$, by
\begin{align}\label{m01}
\begin{split}
f^{(n)} \ = \ \{ ((x_1,\dots,x_n),(y_1,\dots,y_n)) &: (x_1,y_1),\dots, (x_n,y_n) \in f \}, \\
(d^{(n)})((x_1,\dots,x_n),(y_1,\dots,y_n)) \ = \ &\max[d(x_1,y_1), \dots, d(x_n,y_n)].
\end{split}
\end{align}
We will write $(X \times X,d\times d)$ for $(X^{(2)},d^{(2)})$ and $f \times f$ for $f^{(2)}$.

We can regard a chain in $(f^{(n)})^{\times m}$ as an ordered $n$-tuple of chains in $f^{\times m}$.

\begin{lem}\label{lemm01} Let $f$ be a closed relation on a compact metric space $(X,d)$ and $n \geq 2$. Let
$(x_1,\dots,x_n),(y_1, \dots, y_n) \in X^{(n)}$
\begin{itemize}
\item[(a)] If $\A_d f = X \times X$  then
\begin{align}\label{m02a}
\begin{split}
\ell^{f^{(n)}}_{d^{(n)}}((x_1, \dots, x_n),&(x_1, \dots, x_n)) \ = \ 0, \\
\ell^{f^{(n)}}_{d^{(n)}}((y_1, \dots, y_n),&(x_1, \dots, x_n)) \ \leq \\
(n-1) \ell^{f^{(n)}}_{d^{(n)}}((x_1, \dots, &x_n),(y_1, \dots, y_n)).
\end{split}
\end{align}
\item[(b)] If $\CC_d f = X \times X$ then
\begin{align}\label{m02b}
\begin{split}
m^{f^{(n)}}_{d^{(n)}}((x_1, \dots, x_n),&(x_1, \dots, x_n)) \ = \ 0, \\
m^{f^{(n)}}_{d^{(n)}}((y_1, \dots, y_n),&(x_1, \dots, x_n)) \ \leq \\
(n-1) m^{f^{(n)}}_{d^{(n)}}((x_1, \dots, &x_n),(y_1, \dots, y_n)).
\end{split}
\end{align}
\end{itemize}\end{lem}

{\bfseries Proof:} (a): 
%
Let $\ep > 0$. Let $x_{n+1} = x_1$.  Because $\A_d f = X \times X$,
there exist $C_i \in f^{\times m_i}$ such that $||x_i C_i x_{i+1}|| < \ep $ for $i = 1, \dots, n$. The concatenations
$C_i \cdot C_{i+1}\dots C_n \cdot C_1, \dots C_{i-1} \in f^{\times m}$ with $m = \sum_i m_i$ and
\begin{equation}\label{m02c}
||x_i C_i \cdot C_{i+1}\dots C_n \cdot C_1, \dots C_{i-1} x_i||  \ \leq \ n \ep.
\end{equation}

Now let $L = \ell^{f^{(n)}}_{d^{(n)}}((x_1,\dots,x_n),(y_1, \dots, y_n))$ and $x_{n+1} = x_1$.
For $\ep > 0$ there exist $C_i \in f^{\times m}$ such that
$||x_i C_i y_i|| \leq L + \ep$ for $i = 1, \dots, n$. There exist
$D_i \in f^{\times m_i}$ such that $||y_i D_i x_{i+1}|| \leq \ep$ for $i = 1, \dots, n$.
\begin{align}\label{m02d}
\begin{split}
||y_i D_i \cdot C_{i+1} \cdot D_{i+1} \dots D_n \cdot &C_{n} \cdot D_1
\cdot C_1 \dots C_{i-1} \cdot D_{i-1} x_i || \\ &\leq (n-1)L + (2n-1) \ep.
\end{split}
\end{align}

Futhermore, all the concatenations have length $(n-1)m + \sum_i m_i$.

As $\ep > 0$ was arbitrary, (\ref{m02a}) follows.

(b) The estimate (\ref{m02b}) is proved exactly the same way. Simply replace
the chain-length statements using $|| \cdot ||$ with chain-bound statements using $| \cdot |$.

$\Box$ \vspace{.5cm}

\begin{theo}\label{theom01a} Let $f$ be a closed relation on a compact metric space $(X,d)$.
\begin{itemize}
\item[(a)] If $\A_d f = X \times X$ then
$\ell^{f\times f}_{d \times d}$ is a pseudo-metric on $X \times X$. The relation $\A_{d^{(n)}}(f^{(n)})$ is a
closed equivalence relation on $X^{(n)}$ for all $n \in \N$.
\item[(b)] If $\CC_d f = X \times X$ then
$m^{f\times f}_{d \times d}$ is a pseudo-ultrametric on $X \times X$. The relation $\CC_{d^{(n)}}(f^{(n)})$ is a
closed equivalence relation on $X^{(n)}$ for all $n \in \N$.
\end{itemize}\end{theo}

{\bfseries Proof:} With $n = 2$, (\ref{m02a}) and (\ref{m02b}) imply that $\ell^{f\times f}_{d \times d}$ and
$m^{f\times f}_{d \times d}$ are symmetric and vanish on the diagonal.
They satisfy the triangle inequality by  (\ref{c07}).
So, both are is a pseudo-metrics. Since $m^{f\times f}_{d \times d}$
is a pseudo-metric, we have\\ $m^{f\times f}_{d \times d}((y_1,y_2),(y_1,y_2)) = 0$
for all $(y_1,y_2) \in X \times X$. From
(\ref{c07a}) we obtain the ultrametric version of the triangle inequality.

For any $n$, (\ref{m02a}) implies that $((x_1, \dots, x_n),(x_1, \dots, x_n)) \in \A_{d^{(n)}}(f^{(n)})$ and
$((y_1, \dots, y_n),(x_1, \dots, x_n)) \in \A_{d^{(n)}}(f^{(n)})$
if $((x_1, \dots, x_n),(y_1, \dots, y_n)) \in \A_{d^{(n)}}(f^{(n)})$.
Thus, $\A_{d^{(n)}}(f^{(n)})$ is reflexive and symmetric as well as
transitive. Using (\ref{m02b}) we similarly obtain that
$\CC_{d^{(n)}}(f^{(n)})$ is an equivalence relation when $\CC_d f = X \times X$.

$\Box$ \vspace{.5cm}

If $f$ is a closed relation on a compact metric space $(X,d)$ with $\A_d f = X \times X$
we define $\rho_d^f : X \times X \to \R$ by
\begin{equation}\label{m03}
\rho_d^f(x,y) \ = \ \ell^{f \times f}_{d \times d}((z,z),(x,y)) \ = \ \ell^{f \times f}_{d \times d}((x,y),(z,z)),
\end{equation}
Because $\A_d f = X \times X$,  $\ell^{f \times f}_{d \times d}((z,z),(w,w)) = 0$
for all $z,w \in X$. Hence, the above definition is
well-defined, independent of the choice of $z$.

Similarly, if $f$ is a closed relation on a compact metric space $(X,d)$ with $\CC_d f = X \times X$
we define $\theta_d^f : X \times X \to \R$ by
\begin{equation}\label{m04}
\theta_d^f(x,y) \ = \ m^{f \times f}_{d \times d}((z,z),(x,y)) \ = \ m^{f \times f}_{d \times d}((x,y),(z,z)).
\end{equation}
Again the definition is independent of the choice of $z \in X$.

\begin{theo}\label{theom02} Let $f$ be a closed relation on a compact metric space $(X,d)$.
\begin{itemize}
\item[(a)] If $\A_d f = X \times X$ then
$\rho_d^f$ is a pseudo-metric on $X$ with $\rho_d^f \leq d$.
The relation $f$ is an isometry on $(X,\rho_d^f)$.
\item[(b)] If $\CC_d f = X \times X$ then
$\theta_d^f$ is a pseudo-ultrametric on $X$ with $\theta_d^f \leq d$.
The relation $f$ is an isometry on $(X,\theta_d^f)$.
\end{itemize}
\end{theo}

{\bfseries Proof:} (a) By  (\ref{c06})
\begin{equation}\label{m05}
\ell^{f \times f}_{d \times d}((z,z),(x,y)) \leq
\ell^{f \times f}_{d \times d}((z,z),(x,x)) + (d \times d)((x,x),(x,y)).
\end{equation}
Since $\ell^{f \times f}_{d \times d}((z,z),(x,x)) = 0$ and
$(d \times d)((x,x),(x,y)) = d(x,y)$ it follows that $\rho^f_d \leq d$.

Since $\ell^{f \times f}_{d \times d}((z,z),(x,x)) = 0$, $\rho^f_d(x,x) = 0$ for all $x$.

The twist map $(x,y) \mapsto (y,x)$ is a $d \times d$ isometry which maps $f \times f$ to itself. It follows that
$\ell^{f \times f}_{d \times d}((z,z),(x,y)) = \ell^{f \times f}_{d \times d}((z,z),(y,x))$. That is,
$\rho^f_d$ is symmetric.

Now let $L_1 = \rho^f_d(x,y), L_2 = \rho^f_d(y,z)$ and let $\ep > 0$. There
exist $C_{yy}, C_{yx} \in f^{\times n_x}$ such that
$||y C_{yy}y||, ||yC_{yx}x|| < L_1 + \ep$, and $D_{zz}, D_{zy} \in f^{\times n_z}$
such that $||z D_{zz} z||, ||z D_{zy} y|| < L_2 + \ep$.
Finally, there exist $E \in f^{\times n}$ such that $||y E z|| < \ep$. The concatenations
$E \cdot D_{zy} \cdot C_{yx}$ and $C_{yy} \cdot E \cdot D_{zz}$ are in $f^{\times (n + n_x + n_z)}$ and
\begin{equation}\label{m06}
||y E \cdot D_{zy} \cdot C_{yx}x||, ||y C_{yy} \cdot E \cdot D_{zz}z|| < L_1 + L_2 + 3 \ep.
\end{equation}
It follows that $\rho^f_d(x,z) = \ell^{f \times f}_{d \times d}((y,y),(x,z)) \leq L_1 + L_2$.
That is, $\rho^f_d$ satisfies the triangle
inequality.

If $(x_1,y_1) \in \A_{d \times d}^{f \times f}(x,y)$ then
$\ell^{f \times f}_{d \times d}((x,y),(x_1,y_1)) = 0$ and
so \\ $\ell^{f \times f}_{d \times d}((z,z),(x,y)) =
\ell^{f \times f}_{d \times d}((z,z),(x_1,y_1))$ by the triangle inequality.
In particular, if $(x,x_1), (y,y_1) \in f$ then $((x,y),(x_1,y_1))
\in f \times f \subset \A_{d \times d}^{f \times f}$ and so
$\rho^f_d(x,y) = \rho^f_d(x_1,y_1)$.

(b) The proof for the ultrametric version of the triangle inequality
is the only part which requires some adjustment. Let
$L_1 = \theta^f_d(x,y), L_2 = \theta^f_d(y,z)$ and let $\ep > 0$.
Choose $C_{yy},C_{yx}, D_{zz}, D_{zy}, E$ chains as above with
the chain-length inequalities replaced by chain-bound inequalities.
In addition, there exists $\bar E \in f^{\times m}$
such that $|y \bar E y| < \ep$. The concatenations
$E \cdot D_{zy} \cdot \bar E \cdot C_{yx}, \ C_{yy} \cdot \bar E
\cdot E \cdot D_{zz} \in f^{\times (m + n + n_x + n_z)}$ and
as in (\ref{c02a})
\begin{equation}\label{m07}
|y E \cdot D_{zy} \cdot \bar E \cdot C_{yx}x|, |y C_{yy} \cdot
\bar E \cdot E \cdot D_{zz}z| < \max[L_1 + 2 \ep,L_2 + 2 \ep].
\end{equation}
It follows that $\theta_d^f(x,z) \leq \max(L_1,L_2)$. That is, $\theta^f_d$ satisfies the ultrametric triangle inequality.

As above, if $(x_1,y_1) \in \CC_{d \times d}^{f \times f}(x,y)$ then
$\theta^f_d(x,y) = \theta^f_d(x_1,y_1)$. In particular, this holds if
$(x,x_1), (y,y_1) \in f$.

$\Box$ \vspace{.5cm}

\begin{lem}\label{lemm02a} Let $f$ be a closed relation on a compact metric space $(X,d)$ and let $x,y \in X$.
\begin{itemize}
\item[(a)] If $\A_d f = X \times X$ then there exists $x_1 \in \RR f(x)$ such that
$\rho_d^f(x_1,y) = 0$.
\item[(b)] If $\CC_d f = X \times X$ then there exists $x_1 \in \RR f(x)$ such that
$\theta_d^f(x_1,y) = 0$.
\end{itemize}
\end{lem}

{\bfseries Proof:} (a): For each $i \in \N$ there exists $C_i \in f^{\times n_i}$ such that $||x C_i y|| < 1/i$. Let
$z_i \in f^{n_i}(x)$ and so there exists $D_i \in f^{\times n_i}$ such that $||x D_i z_i || = 0$. Choose a subsequence
$\{ z_{i'} \}$ which converges to $x_1$. Since $||x D_{i'} x_1|| \to 0$, it follows that $\rho_d^f(x_1,y) = 0$.

(b): As before,  choosing $C_i$ so that $|x C_i y| < 1/i$.

$\Box$ \vspace{.5cm}

\begin{lem}\label{lemm03} For a compact metrizable spaces
$X_1, X_2 $, if $\tilde d$ is an admissible metric on $X_1 \times X_2$, then there
exist  admissible metrics $d_1$ on $X_1$ and $d_2$ on $X_2$ such that $\tilde d \leq d_1 \times d_2$
where $(d_1 \times d_2)((x_1,y_1),(x_2,y_2)) = \max[d_1(x_1,x_2), d_2(y_1,y_2)]$ \end{lem}

{\bfseries Proof:} With $A$ a compact metrizable space, on the space of continuous functions $C(A, X_1 \times X_2)$
we let $\tilde d$ denote the sup metric induced by
$\tilde d$ on $X_1 \times X_2$. That is, for
$h_1,h_2 \in C(A, X \times X)$, $\tilde d(h_1,h_2) = \sup \{ \tilde d(h_1(a),h_2(a)) : a \in A \}$.
The maps $i : X_1 \to C(X_2,X_1 \times X_2)$ and
$j : X_2 \to C(X_1,X_1 \times X_2)$ defined by $i(x)(y) = (x,y), j(y)(x) = (x,y)$ are continuous
injections. Define $d$ on $X$ by
\begin{equation}\label{m07a}
d_1(x_1,x_2) = 2 \tilde d(i(x_1),i(x_2)), \quad d_2(y_1,y_2) = 2 \tilde d(j(y_1),j(y_2))].
\end{equation}
 If $(x_1,y_1),(x_2,y_2) \in X_1 \times X_2$
then
\begin{equation}\label{m08}
\begin{split}
\tilde d((x_1,y_1),(x_2,y_2)) \ \leq \ \tilde d((x_1,y_1),(x_2,y_1)) + \tilde d((x_2,y_1),(x_2,y_2)) \ \leq \\
\tilde d(i(x_1),i(x_2)) + \tilde d(j(y_1),j(y_2)) \ \leq \ \max[d_1(x_1,x_2),d_2(y_1,y_2)].
\end{split}
\end{equation}

$\Box$ \vspace{.5cm}

\begin{theo}\label{theom04} If $f$ is a closed relation on a compact metrizable space, then
$\G (f^{(n)}) \ = \ \bigcap_d  \ \{ \A_{d^{(n)}} (f^{(n)}) \}$ with
$d$ varying over the set of admissible metrics on $X$.
Furthermore, there exists an admissible metric $d$ on $X$ such that
$\G (f^{(n)}) \ = \ \A_{d^{(n)}} (f^{(n)}) $ for all $n \in \N$. \end{theo}

{\bfseries Proof:} By Theorem \ref{theoc01}
$\G (f^{(n)}) \ = \ \bigcap_{\tilde d}  \ \{ \A_{\tilde d} (f^{(n)}) \}$ as $\tilde d$
varies over the admissible metrics on $X^{(n)}$. For any such
 $\tilde d$, Lemma \ref{lemm03} implies there exist admissible metrics
$d_1, \dots, d_n$ on $X$
such that $d_1 \times \dots d_n \geq \tilde d$. With $d = \max_i d_i$, $d^{(n)} \geq \tilde d$.
Then $\A_{d^{(n)}} (f^{(n)}) \subset \A_{\tilde d} (f^{(n)})$ and so we need only
intersect over the metrics of the form $d^{(n)}$ with $d$ an admissible metric on $X$.

Thus, $\G (f^{(n)}) \ \subset \ \A_{d^{(n)}} (f^{(n)}) $ and if
$((x_1,\dots,x_n),(y_2, \dots y_n)) \not\in \G (f^{(n)})$ then there exists $d$ so that
$\ell_{d^{(n)}}^{f^{(n)}}((x_1,\dots,x_n),(y_2, \dots y_n)) > 0$. Thus,
$\{ \{ \ell_{d^{(n)}}^{f^{(n)}} > 0 \} : d $ an admissible metric on $X \}$ is an
open cover of $[X^{(n)} \times X^{(n)}] \setminus \G(f^{(n)})$. Since
this space is Lindel\"{o}f we can choose a countable subcover, indexed
by a sequence of metrics $\{ d_i : i \in \N \}$ on $X$.  If $d_i$ is
bounded by $M_i$ then $d = \sum_i  \frac{1}{M_i 2^i}d_i$ is an admissible metric on
$X$ and $\A_{d^{(n)}} (f^{(n)}) \subset \bigcap_i \A_{d_i^{(n)}} (f^{(n)}) = \G (f^{(n)})$.

We thus have for each $n \in \N$ a metric $d_n$ such that
$\A_{d_n^{(n)}} (f^{(n)}) = \G (f^{(n)})$. If $d = \sum_n  \frac{1}{M_n 2^n}d_n$
where $M_n$ is a bound on $d_n$, then $d$ is an admissible metric and $\A_{d^{(n)}} (f^{(n)}) = \G (f^{(n)})$ for all $n$.

$\Box$ \vspace{.5cm}

\begin{cor}\label{corm05}  Let $f$ be a closed relation on a compact metric space $(X,d)$.
 If $\G f = X \times X$ then
the relation $\G(f^{(n)})$ is a
closed equivalence relation on $X^{(n)}$ for all $n \in \N$. \end{cor}

{\bfseries Proof:} Let $d$ be a metric on $X$ such that $\G (f^{(n)}) \ = \ \A_{d^{(n)}} (f^{(n)}) $.
Since $X \times X = \G f \subset \A_d f$, Theorem \ref{theom01a}
implies that $\G(f^{(n)})$ is a closed equivalence relation.

$\Box$ \vspace{.5cm}

We define the \emph{synchrony relations}.

If $f$ is a closed relation on a compact metric space $(X,d)$ with $\A_d f = X \times X$ we define the relation on $X$:
\begin{equation}\label{m09}
R_{\ell}f \ = \ \A_{d \times d}(f \times f)(1_X) \ =
\ \A_{d \times d}(f \times f)^{-1}(1_X) \ = \ \{ \rho^f_d = 0 \}.
\end{equation}

If $\CC_d f = X \times X$ we define the  relation on $X$:
\begin{equation}\label{m10}
R_{m}f \ = \ \CC_{d \times d}(f \times f)(1_X) \ =
\ \CC_{d \times d}(f \times f)^{-1}(1_X) \ = \ \{ \theta^f_d = 0 \}.
\end{equation}

As these are zero-sets of a continuous pseudo-metrics, they are closed equivalence relations.

If $\G f = X \times X$ we define the  relation on $X$:
\begin{equation}\label{m11}
R_{g}f \ = \ \G (f \times f)(1_X) \ = \ \G (f \times f)^{-1}(1_X).
\end{equation}

Since $\G (f \times f) = \A_{d \times d}(f \times f)$ for some admissible
metric $d$, $R_{g} f$ is a closed equivalence relation as well.

\begin{prop}\label{propm05a} If $f$ is a closed relation on a compact metric space
$(X,d)$ with $\A_d f = X \times X$ then
$\A_d R_{\ell} = R_{\ell}$. If  $\CC_d f = X \times X$ then
$\CC_d R_{m} = R_{m}$. \end{prop}

{\bfseries Proof:} Suppose $C = [(a_1,b_1),\dots, (a_n,b_n)] \in (R_{\ell})^n$ and
$||x C y|| < \ep$. Let $b_0 = x, a_{n+1} = y$.
 For $i = 1, \dots , n$ $\rho_d^f(a_i,b_i) = 0$. For $i = 1,\dots,n+1$,
 $\rho_d^f(b_i,a_{i-1}) \leq d(b_i,a_{i-1}).$
 It follows from the triangle inequality for $\rho_d^f$ that
 $\rho_d^f(x,y) \leq ||x C y|| < \ep$. As $\ep$ was arbitrary, it follows
 that $(x,y) \in R_{\ell}$.

 Since $\theta_d^f$ is a pseudo-ultrametric, it similarly follows that if
 $C \in (R_m)^n$ then $\theta_d^f(x,y) \leq |x C y|$.

$\Box$ \vspace{.5cm}

\begin{lem}\label{lemm06} Let $(x_1,\dots,x_n) \in X^{(n)}$, $y \in X$.
\begin{enumerate}
\item[(a)]  If $\A_d f = X \times X$ and  $x_1,\dots, x_n \in R_{\ell} f(y)$,
then\\ $\ell_{d^{(n)}}^{f^{(n)}}((y,\dots,y), (x_1,\dots,x_n)) = 0$.
\item[(b)]  If $\CC_d f = X \times X$ and  $x_1,\dots, x_n \in R_{m} f(y)$,
then\\ $m_{d^{(n)}}^{f^{(n)}}((y,\dots,y), (x_1,\dots,x_n)) = 0$.
\end{enumerate}
\end{lem}

{\bfseries Proof:} (a): Given $\ep > 0$, there exist $C_i, D_i \in f^{\times m_i}$
such that $||y C_i x_i|| < \ep$ and $||y D_i y|| < \ep$.
\begin{equation}\label{m12}
||y D_1 \cdot D_2 \dots \widehat{D_i}\dots D_n \cdot C_i x_i|| < n \ep,
\end{equation}
where $\widehat{D_i}$ denotes the omission of $D_i$ in the sequence.

These concatenations have length $\sum_i m_i$. So the result follows.

(b): As usual replace the chain-length estimates by chain-bound estimates.

$\Box$ \vspace{.5cm}

\begin{prop}\label{propm07}  Let $f$ be a closed relation on a compact
metric space $(X,d)$ and let $x,y \in X$, $N \in \N$ and $\ep > 0$.
\begin{enumerate}
\item[(a)]  If $\A_d f = X \times X$ there exists $n > N$ such that
for all $x' \in R_{\ell}(x), y' \in R_{\ell}(y)$ there exists $C \in f^{\times n}$
such that $|| x' C y' || < \ep$.
\item[(b)]  If $\CC_d f = X \times X$ there exists $n > N$ such that
for all $x' \in R_m(x), y' \in R_m(y)$ there exists $C \in f^{\times n}$
such that $| x' C y' | < \ep$.
\end{enumerate}
\end{prop}

{\bfseries Proof:} (a): Choose $\{ x_1, \dots, x_{m_1} \} $ an $\ep/6$ dense subset
of $ R_{\ell} f(x)$, i.e. every point of  $ R_{\ell} f(x)$
has distance less than $\ep/6$ from some $x_i$ and choose $\{ y_1, \dots, y_{m_2} \} $
an $\ep/6$ dense subset of $ R_{\ell} f(y)$.
By Lemma \ref{lemm06} $\ell_{d^{(m_1)}}^{f^{(m_1)}}((x_1,\dots,x_{m_1}),(x,\dots,x)) = 0$ and
$\ell_{d^{(m_2)}}^{f^{(m_2)}}((y, \dots,y),(y_1,\dots,y_{m_2})) = 0$. Also,
$\ell_d^f(x,x) = \ell_d^f(x,y) = 0$. Choose
$C_i \in f^{\times k_1}, D_j \in f^{\times k_2}$ such that $||x_i C_i x|| < \ep/6$
and $||y D_j y_j|| < \ep/6$ for $i = 1,\dots,m_1, j= 1,\dots, m_2$.
There exist $E \in f^{\times k_3}, F \in f^{\times k_4}$ such that
$||x E x|| < \ep/6N$ and $|| x F y|| < \ep/6$.
 \begin{equation}\label{m13a}
 || x_i C_i \cdot E^N \cdot F \cdot D_j y_j|| \ < \ 4\ep/6.
 \end{equation}
 The length of each concatenation is $n = k_1 + k_2 + N k_3 + k_4 > N$. If $d(x',x_i) < \ep/6, d(y',y_j) < \ep/6$ then
 $ || x' C_i \cdot E^N \cdot F \cdot D_j y'||  < \ep$.

 (b): Replace the chain-length estimates by chain-bound estimates.

$\Box$ \vspace{.5cm}

\begin{df}\label{defm08} Let $f$ be a closed relation on a compact
metric space $(X,d)$. We call $f$ \emph{chain transitive} when
$\CC_d f = X \times X$, \emph{strong chain transitive} when
$\A_d f = X \times X$, and \emph{vague transitive} when $\G f = X \times X$.
We call $f$ \emph{chain mixing} when $f \times f$ is chain
transitive, \emph{strong chain mixing} when $f \times f$ is strong chain
transitive, and \emph{vague mixing} when $f \times f$ is vague transitive.
\end{df}

Since each projection $\pi : X \times X \to X$ maps $f \times f$,
$\G (f \times f)$, $\A_d (f \times f)$ and $\CC_d (f \times f)$ to
$f$, $\G f$, $\A_d f$ and $\CC_d f$, respectively, it follows
that each sort of mixing implies the corresponding transitivity property.

As observed above, if any of these conditions hold then $f$ is a
surjective relation, by \cite{AW17} Proposition 6.10, 6.11.

\begin{theo}\label{theom09}  Let $f$ be a closed relation on a compact metric space $(X,d)$.
\begin{enumerate}
\item[(a)] If $f$ is chain transitive, then following are equivalent.
\begin{itemize}
\item[(i)] The relation $f$ is chain mixing.
\item[(ii)] For every $n \in \N$, $f^{(n)}$ is chain mixing.
\item[(iii)] $R_m = X \times X$.
\item[(iv)] $f \subset R_m$.
\item[(v)] There exists $x \in X$ such that $(x,y) \in R_m$ for all $y \in f(x)$.
\item[(vi)] For every $\ep > 0$ and $x, y \in X$ there exists $N \in \N$ such that
for all $n \geq N$ there exists $C \in f^{\times n}$ such that $|x C y| < \ep$.
\item[(vii)] For every $\ep > 0$ there exists $N \in \N$ such that
for all $x, y \in X$ and all $n \geq N$ there exists $C \in f^{\times n}$ such that $|x C y| < \ep$.
\end{itemize}

\item[(b)] If $f$ is strong chain transitive, then following are equivalent.
\begin{itemize}
\item[(i)] The relation $f$ is strong chain mixing.
\item[(ii)] For every $n \in \N$, $f^{(n)}$ is strong chain mixing.
\item[(iii)] $R_{\ell} = X \times X$.
\item[(iv)] $f \subset R_{\ell}$.
\item[(v)] There exists $x \in X$ such that $(x,y) \in R_{\ell}$ for all $y \in f(x)$.
\item[(vi)] For every $\ep > 0$ and $x, y \in X$ there exists $N \in \N$ such that
for all $n > N$ there exists $C \in f^{\times n}$ such that $||x C y|| < \ep$.
\item[(vii)] For every $\ep > 0$ there exists $N \in \N$ such that
for all $x, y \in X$ and all $n > N$ there exists $C \in f^{\times n}$ such that\\ $||x C y|| < \ep$.
\end{itemize}

\item[(c)] If $f$ is vague transitive, then following are equivalent.
\begin{itemize}
\item[(i)] The relation $f$ is vague mixing.
\item[(ii)] For every $n \in \N$, $f^{(n)}$ is vague mixing.
\item[(iii)] $R_g = X \times X$.
\item[(iv)] $f \subset R_g$.
\item[(v)] There exists $x \in X$ such that $(x,y) \in R_g$ for all $y \in f(x)$.
\end{itemize}
\end{enumerate}
\end{theo}

{\bfseries Proof:} We prove (b). As usual, the proofs of the equivalences of (a) are completely
analogous.

(ii) $\Rightarrow$ (i) : Obvious.

(i) $\Rightarrow$ (iii) : $\A_{d \times d} (f \times f) = X^{(2)} \times X^{(2)}$ says that
for all $(x,y), (z,w) \in X^{(2)}$ $\ell_{d \times d}^{f \times f}((z,w),(x,y)) = 0$. Hence,
for all $x,y \in X$, $\rho_d^f(x,y) = 0$.

(iii) $\Rightarrow$ (iv) $\Rightarrow$ (v) : Obvious.

(v) $\Rightarrow$ (iii) : By  Theorem \ref{theom02} $f$ is an isometry on $(X,\rho_d^f)$.
If $y_n \in f^{n}(x)$ then there exists $y_1, \dots, y_n$ such that $y_1 \in f(x)$ and
$y_i \in f(y_{i-1})$ for $i = 2, \dots, n$. Hence, by induction, $\rho_d^f(y_i,y_{i-1}) = \rho_d^f(x,y_1) = 0$
for $i = 2, \dots, n$. Hence, $\rho_d^f(x,y_i) = 0$ for $i = 1,\dots, n$. It follows that
$\bigcup_n \ f^n(x) \subset R_{\ell}(x)$.  Since $R_{\ell}(x)$ is closed, $\RR f(x) \subset R_{\ell}(x)$.
For any $y_1, y_2 \in X$ there exist $x_1,x_2 \in \RR f(x)$ such that $(x_1,y_1), (x_2,y_2) \in R_{\ell}$.
Since  $\RR f(x) \subset R_{\ell}(x)$, $(x,x_1), (x,x_2) \in R_{\ell}$. Since $R_{\ell}$ is an
equivalence relation $(y_1,y_2) \in R_{\ell}$.

(iii) $\Rightarrow$ (vii) : Since $ R_{\ell}(x) = R_{\ell}(y) = X$, Proposition \ref{propm07} implies
there exists $N \in \N$ so that for all $x', y$ there exists $C \in f^{\times N}$ such that
$||x' C y|| < \ep$. Given $n > N$, choose $x' \in f^{n-N}(x)$. Thus, there exists $D \in f^{\times (n-N)}$
such that $||x D x' || = 0$. Hence, $D \cdot C \in f^{\times n}$ with $|| x  D \cdot C y || < \ep$.

(vii) $\Rightarrow$ (vi) : Obvious.

(vi) $\Rightarrow$ (ii) : Given $\ep > 0$ and $(x_1,\dots,x_n), (y_1,\dots,y_n) \in X^{(n)}$,
there exist $N_i \in \N$
such that if $m > N_i$ there exists $C_i \in f^{\times m}$ such that $||x_i C_i y_i|| < \ep$. Choose
$m > \max_i N_i$. As $\ep > 0$ was arbitrary, it follows that
$\ell_{d^{(n)}}^{f^{(n)}}((x_1,\dots,x_n), (y_1,\dots,y_n)) = 0$. It follows that $\A_{d^{(n)}}(f^{(n)}) =
X^{(n)} \times X^{(n)}$ for all $n$. Replacing $n$ by $2n$ we see that $f^{(n)}$ is strong chain mixing.

(c): Choose an admissible metric $d$ on $X$ such that $\G (f \times f) = \A_{d \times d}(f \times f)$
and apply (b).

$\Box$ \vspace{.5cm}

Now we consider when the mixing conditions fail.

\begin{theo}\label{theom10} Let $f$ be a closed relation on a compact metric space $(X,d)$ and
let $g$ be a continuous map on a metrizable space $Y$. Assume that $\pi : X \to Y$ is a continuous
surjection mapping $f$ to $g$.
\begin{enumerate}
\item[(a)] Assume that $\CC_d f = X \times X$, that $Y$ is zero-dimensional and that $g$ is
equicontinuous. Then $g$ is a minimal homeomorphism on $Y$ and so $(Y,g)$ is a single
periodic orbit if $Y$ is finite or an odometer if $Y$ is infinite. Furthermore,
$R_m \subset (\pi \times \pi)^{-1}(1_Y)$.

\item[(b)] Assume that $\A_d f = X \times X$, that $g$ is an isometry for $(Y,\bar d)$ with $\bar d$ an
admissible metric on $Y$ and that
$\pi : (X,d) \to (Y,\bar d)$ is Lipschitz. Then $g$ is a minimal homeomorphism on
$Y$ and $R_{\ell} \subset (\pi \times \pi)^{-1}(1_Y)$.

\item[(c)] Assume that $\G f = X \times X$, and that $g$ is equicontinuous. Then $g$ is a minimal homeomorphism on
$Y$ and $R_{g} \subset (\pi \times \pi)^{-1}(1_Y)$.
\end{enumerate}
\end{theo}

{\bfseries Proof:} In each of these cases $f$ is a surjective relation. Since $\pi$ is surjective, $Y$ is compact.
Since it maps $f$ to $g$, $g$ is a surjective map. In each case the equicontinuity assumption then implies that
$g$ is a homeomorphism by Theorem \ref{theoc03}.

(a): Because $Y$ is zero-dimensional it can be embedded as a closed subset of $\{0,1\}^{\N}$ and so admits
an ultrametric $\bar d$.  By Theorem \ref{theoc03} $g$ is an isometry of $\bar d_f$ which is also a continuous
ultrametric. So replacing $\bar d$ by $\bar d_f$ we can assume that $g$ is an isometric homeomorphism on
the compact ultrametric space $(Y,\bar d)$.  Since $\pi$ maps $f$ to $g$ it maps $\CC_d f$ to
$\CC_{\bar d} g$. Hence, $\CC_{\bar d} g = Y \times Y$. By Theorem \ref{theoc03} again,
$Y = \CC_{\bar d} g(y) = \RR g(y)$ for all $y \in Y$. Hence, $g$ is a minimal homeomorphism.
Since $\pi \times \pi$ maps $f \times f$ to $g \times g$, it maps $\CC_{d \times d} (f \times f)$ to
$\CC_{\bar d \times \bar d} (g \times g)$. For  $x \in X$ let $y = \pi(x)$.
\begin{align}\label{m13}
\begin{split}
(\pi \times \pi)(\CC_{d \times d} (f \times f)(x,x)) \ \subset \ &\CC_{\bar d \times \bar d} (g \times g)(y,y)\\
= \ \RR (g  \times g)(y,y) \ \subset \ &1_Y.
\end{split}
\end{align}
From (\ref{m10}) we see that $R_m f  \subset (\pi \times \pi)^{-1}(1_Y)$.

Because $g$ is a minimal equicontinuous homeomorphism on a zero-dimensional space, it follows that $g$ is
either a periodic orbit or an odometer.  See, e.g. \ \cite{A96} Theorem 3.5.

(b): The proof is similar to that of (a). Since $\pi$ is Lipschitz and maps $f$ to $g$ it maps $\A_d f$ to
$\A_{\bar d} g$. Hence, $\A_{\bar d} g = Y \times Y$. By Theorem \ref{theoc03} again,
$Y = \A_{\bar d} g(y) = \RR g(y)$ for all $y \in Y$. Hence, $g$ is a minimal homeomorphism.
Since $\pi \times \pi$ maps $f \times f$ to $g \times g$, it maps $\A_{d \times d} (f \times f)$ to
$\A_{\bar d \times \bar d} (g \times g)$. For  $x \in X$ let $y = \pi(x)$.
\begin{align}\label{m14}
\begin{split}
(\pi \times \pi)(\A_{d \times d} (f \times f)(x,x)) \ \subset \ &\A_{\bar d \times \bar d} (g \times g)(y,y)\\
= \ \RR (g  \times g)(y,y) \ \subset \ &1_Y.
\end{split}
\end{align}
From (\ref{m09}) we see that $R_{\ell}  f  \subset  (\pi \times \pi)^{-1}(1_Y)$.

(c): Choose $d_1$ a continuous metric on $X$ so that $\G (f \times f) = \A_{d_1 \times d_1}(f \times f)$.
We can choose an admissible metric $\bar d$ on $Y$ so that $g$ is an isometry for $(Y,\bar d)$. Replace
$d$ by the admissible metric defined by $d(x_1,x_2) = d_1(x_1,x_2) + \bar d(\pi(x_1),\pi(x_2))$.
Then $\pi : (X,d) \to (Y, \bar d)$ is Lipschitz. Since
$d \geq d_1$, $\G (f \times f) = \A_{d \times d}(f \times f)$.
 Apply (b) to get  $R_g = R_{\ell} \subset (\pi \times \pi)^{-1}(1_Y)$.

$\Box$ \vspace{.5cm}

\begin{theo}\label{theom11} Let $f$ be a closed relation on a compact metric space $(X,d)$.
\begin{enumerate}
\item[(a)] Assume that $\CC_d f = X \times X$. Let $\pi : (X, \theta_d^f) \to (X/R_m,\bar \theta)$ be the projection to
the quotient metric space and let $ f_m = (\pi \times \pi)(f)$. The quotient space
$(X/R_m,\bar \theta)$ is a compact, ultrametric
space, with $f_m$ a minimal isometric homeomorphism and $\pi: (X,d) \to (X/R_m,\bar \theta)$
is Lipschitz with Lipschitz constant $1$.

If $g$ is an equicontinuous map on a compact, zero-dimensional space $Y$ and
$h : X \to Y$ is continuous mapping $f$ to $g$, then there
exists a continuous $q : X/R_m \to Y$ such that $h = q \circ \pi$ and so $q$ maps $f_m$ to $g$.

\item[(b)] Assume that $\A_d f = X \times X$. Let $\pi : (X, \rho_d^f) \to (X/R_{\ell},\bar \rho)$ be the projection to
the quotient metric space and let $ f_{\ell} = (\pi \times \pi)(f)$. The
quotient space $(X/R_{\ell},\bar \rho)$ is a compact, metric
space, with $f_{\ell}$ a minimal isometric homeomorphism and
$\pi: (X,d) \to (X/R_{\ell},\bar \rho)$ is Lipschitz with Lipschitz constant $1$.

If $g$ is an isometry on a compact, metric space $(Y,d_1)$ and
$h : (X,d) \to (Y,d_1)$ is Lipschitz mapping $f$ to $g$, then there
exists a continuous $q : X/R_{\ell} \to Y$ such that $h = q \circ \pi$ and so $q$ maps $f_{\ell}$ to $g$.

\item[(c)] Assume that $\G f = X \times X$. Let $\pi : X \to X/R_{g}$ be the projection to
the space of $R_g$ equivalence classes, and let $ f_{g} = (\pi \times \pi)(f)$.
The quotient space $X/R_{g}$ is a compact, metrizable
space, with $f_{g}$ a minimal equicontinuous homeomorphism.

If $g$ is an equicontinuous map on a compact, metrizable space $Y$ and $h : X \to Y$
is continuous mapping $f$ to $g$, then there
exists a continuous $q : X/R_g \to Y$ such that $h = q \circ \pi$ and so $q$ maps $f_{g}$ to $g$.
\end{enumerate}
\end{theo}

{\bfseries Proof:} (a): The relation $f$ is an isometry of $(X,\theta_d^f)$ and
$\theta_d^f$ is a pseudo-ultrametric with $d \geq \theta_d^f$
by  Theorem \ref{theom02}. By Proposition \ref{propc04} and Theorem \ref{theoc03}
the induced relation $f_m$ is an isometric, minimal homeomorphism
on the compact ultrametric space $(X/R_m,\bar \theta)$.

Because $g$ is a map, $(h \times h)(f)$ is the restriction of $g$ to the
$g$ $^+$invariant set $h(X)$ and so, by replacing $Y$ by $h(X)$ and
$g$ by its restriction, we may assume that $h$ is surjective. By
Theorem \ref{theom10} (a) $R_m \subset (h \times h)^{-1}(1_Y)$.
Thus, $h$ is constant on each $R_m$ equivalence class and so factors to a continuous map $q : X/R_m \to Y$. Furthermore,
\begin{equation}\label{m15}
(q \times q)(f_m) \ = \ (q \times q)(\pi \times \pi)(f) \ = \ (h \times h)(f) \ = \ g.
\end{equation}
That is, $q$ maps $f_m$ to $g$.

(b): The relation $f$ is an isometry of $(X,\rho_d^f)$ and $\rho_d^f$ is a pseudo-metric with $d \geq \theta_d^f$
by  Theorem \ref{theom02}. By Proposition \ref{propc04} and Theorem
\ref{theoc03} the induced relation $f_{\ell}$ is an isometric, minimal homeomorphism
on the compact metric space $(X/R_{\ell},\bar \rho)$.

Again, by replacing $Y$ by $h(X)$ and
$g$ by its restriction, we may assume that $h$ is surjective. By
Theorem \ref{theom10} (b) $R_{\ell} \subset (h \times h)^{-1}(1_Y)$.
Thus, $h$ is constant on each $R_{\ell}$ equivalence class and so factors
to a continuous map $q : X/R_{\ell}\to Y$. As above
$q$ maps $f_{\ell}$ to $g$.

(c): As usual, choose a continuous metric $d$ so that
$\G (f \times f) = \A_{d \times d}(f \times f)$. Then $R_g = R_{\ell}$
and so $f_g = f_{\ell}$ is a minimal isometric homeomorphism of the compact metric space $(X/R_g, \bar \theta)$ by (b).

Again we may assume that $h$ is surjective and obtain
$R_{g} \subset (h \times h)^{-1}(1_Y)$ from Theorem \ref{theom10} (c). Complete
the proof as before.

$\Box$ \vspace{.5cm}

For $A \subset X$ and $d$ a metric on $X$, let
$d(A,x) = \inf \{ d(x_1,x) : x_1 \in A \}$. If $A$ is closed, then compactness implies that
there exists $x_1 \in A$ such that $d(x_1,y) = d(A,y)$. In particular,  $d(A,y) = 0$ iff $y \in A$.

\begin{lem}\label{lemm13} If a map $f$ is an isometry on a compact metric
space $(X,d)$, then $\ell_d^f(x,y) = d(\RR f(x),y)$ and
if $d$ is an ultrametric then $m_d^f(x,y) = d(\RR f(x),y)$.  In particular, $\A_d f = \RR f$ and if $d$ is an
ultrametric $\CC_d f = \RR f$. \end{lem}

{\bfseries Proof:} In general, if $x_1 \in \RR f(x)$ and $\ep > 0$, then $\ell_d^f(x,x_1) = m_d^f(x,x_1) = 0$ since
$\RR f \subset \A_d f \subset \CC_d f$. Hence, by (\ref{c06}) $\ell_d^f(x,y), m_d^f(x,y) \leq d(x_1,y)$.

Let $L = \ell_d^f(x,y)$. Given $\ep > 0$ there
exists $C \in f^{\times n}$ such that $||x C y|| < L + \ep$. Let
$C = [(a_1,f(a_1)), \dots, (a_n,f(a_n))$ and let $a_{n+1} = y$.
Let $d(x,a_1) = \ep_0,d(f(a_i),a_{i+1}) = \ep_i$ for $i = 1, \dots, n$ so
that $\sum_i \ep_i < L + \ep$. Since $f$ is an isometry,
$d(f^n(x),f^n(a_1)) = \ep_0,d(f^{n+1-i}(a_i),f^{n-i}a_{i+1}) = \ep_i$ for $i = 1, \dots, n$. By the triangle inequality,
 $d(f^n(x),a_{n+1}) = d(f^n(x),y) \leq \sum_i \ep_i < L + \ep$. That is, for every $k \in \N$ there exists $n_k$ such that
 $d(f^{n_k}(x),y) < L + (1/k)$. By going to a subsequence we can assume
 that $\{ f^{n_k}(x) \}$ converges to a point $x_1 \in \RR f(x)$
 and so $d(x_1,y) \leq L$.

 Now assume that $d$ is an ultrametric and $L = m_d^f(x,y)$. Given $\ep > 0$ there
exists $C \in f^{\times n}$ such that $|x C y| < L + \ep$. Proceeding
as above we have that $\max_i \ep_i < L + \ep$. This time the
ultrametric version of the triangle inequality implies that
$d(f^n(x),y) \leq \max_i \ep_i < L + \ep$. Again we can choose a convergent subsequence
to obtain $x_1 \in \RR f(x)$ such that $d(x_1,y) \leq L$.

In particular, if $y \in \A_d f(x)$, then $\ell_d^f(x,y) = 0$ and so $d(\RR f(x),y) = 0$. Hence, $y \in \RR f(x)$.
Similarly, if $d$ is an ultrametric, then
$y \in \CC_d f(x)$ implies $y \in \RR f(x)$. See also Theorem \ref{theoc03} (d).

$\Box$ \vspace{.5cm}

\begin{cor}\label{corm11a} Let $f$ be a closed relation on a compact metric space $(X,d)$.
\begin{enumerate}
\item[(a)] The following are equivalent.
\begin{itemize}
\item[(i)]The relation $f$ is a minimal, equicontinuous homeomorphism with
$X$ zero-dimensional (and so is an odometer or a periodic orbit).
\item[(ii)]$\CC_d f = X \times X$ and $R_{m} = 1_X$.
\item[(iii)] $\CC_d f = X \times X$  and for all $n \in \N$, $\CC_{d^{(n)}} f^{(n)} = \RR f^{(n)}$.
\end{itemize}
\item[(b)]  The following are equivalent.
\begin{itemize}
\item[(i)]The relation $f$ is a minimal homeomorphism with a continuous metric $\rho \leq d$ such that $f$ is
an isometry on $(X,\rho)$.
\item[(ii)]$\A_d f = X \times X$ and $R_{\ell} = 1_X$.
\item[(iii)] $\A_d f = X \times X$  and for all $n \in \N$, $\A_{d^{(n)}} f^{(n)} = \RR f^{(n)}$.
\end{itemize}

 \item[(c)]   The following are equivalent.
\begin{itemize}
\item[(i)] The relation $f$ is a minimal, equicontinuous homeomorphism.
\item[(ii)] $\G f = X \times X$ and $R_{g} = 1_X$.
\item[(iii)] $\G f = X \times X$  and for all $n \in \N$, $\G f^{(n)} = \RR f^{(n)}$.
\end{itemize}
\end{enumerate}
\end{cor}

{\bfseries Proof:} (a) (ii) $\Rightarrow$ (i): If $\CC_d f = X \times X, \A_d f = X \times X,$ or $\G f = X \times X$ and
the corresponding synchrony relation is $1_X$
then $f$ is a minimal, equicontinuous homeomorphism with the conditions in (i) following
from the first parts of each section of Theorem \ref{theom11} since the quotient map $\pi$ is the identity and
$f = f_m, f_{\ell}$ or $f_g$, respectively.

(i)$\Rightarrow$ (ii): If $f$ satisfies the equicontinuity assumption, then
minimality implies $\G f = X \times X$, etc. and
 from the second parts of each section of Theorem \ref{theom11}
it follows that the identity map on $(X,d)$ factors through the quotient map
$\pi : X \to X/R$.   This requires $R \subset 1_X$.

(iii) $\Rightarrow$ (ii): $\CC_d f = X \times X, \A_d f = X \times X,$ or $\G f = X \times X$ implies
$\CC_{d \times d} (f \times f)(1_X) \supset 1_X$, etc.  Since
$\RR (f \times f)(1_X) \subset 1_X$, it follows that $R = 1_X$ in each case.

(i)$\Rightarrow$ (iii): (a) Since $f$ is minimal, $\CC_d f = X \times X$. Since $d \geq \theta_d^f$, it follows that
$\RR f^{(n)} \subset \CC_{d^{(n)}} f^{(n)} \subset \CC_{\theta^{(n)}} f^{(n)}$.
Because $\theta^{(n)}$ is an ultrametric on
$X^{(n)}$ it follows from Lemma \ref{lemm13} that $\CC_{\theta^{(n)}} f^{(n)} = \RR f^{(n)}$.

(b) Since $f$ is minimal, $\A_d f = X \times X$. Since $d \geq \rho_d^f$, it follows that
$\RR f^{(n)} \subset \A_{d^{(n)}} f^{(n)} \subset \A_{\rho^{(n)}} f^{(n)}$. Again Lemma \ref{lemm13}
implies that $\A_{\theta^{(n)}} f^{(n)} = \RR f^{(n)}$.

 (c) Choose a metric $d$ so that $\G f^{(n)} = \A_{d^{(n)}} f^{(n)}$ for all
 $n \in \N$. Since $f$ is minimal $\G f = \A_d f = X \times X$.
 Since $f$ is equicontinuous, $d_f $ is a continuous metric on $X$ and $f$ is
 an isometry on $(X,d_f)$. Since $d_f \geq d$,
 $\G f^{(n)} \subset \A_{d_f^{(n)}} f^{(n)} \subset \A_{d^{(n)}} f^{(n)} = \G  f^{(n)}$.
 Applying (b) to the isometry $f$ on $(X,d_f)$
 we have that $\G  f^{(n)} = \A_{d_f^{(n)}} f^{(n)} = \RR f^{(n)}$.

$\Box$ \vspace{.5cm}

In addition, we obtain the following dichotomy result.

\begin{cor}\label{corm12} Let $f$ be a closed relation on a compact metric space $(X,d)$.
\begin{enumerate}
\item[(a)] Assume that $f$ is chain transitive. Exactly one of the following is true:
\begin{itemize}
\item[(i)] The relation $f$ is chain mixing.
\item[(ii)] There exists a continuous $\pi : X \to Y$ mapping $f$ to $g$ on $Y$ with $Y$
consisting of a single, nontrivial periodic
orbit of $g$.
\end{itemize}

\item[(b)] Assume that $f$ is strong chain transitive. Exactly one of the following is true:
\begin{itemize}
\item[(i)] The relation $f$ is strong chain mixing.
\item[(ii)] There exists a Lipschitz $\pi : (X,d) \to (Y,d_1)$ mapping $f$ to $g$ on $Y$
with $g$ a minimal isometric homeomorphism
on the nontrivial compact metric space $(Y,d_1)$.
\end{itemize}

\item[(b)] Assume that $f$ is vague transitive. Exactly one of the following is true:
\begin{itemize}
\item[(i)] The relation $f$ is vague mixing.
\item[(ii)] There exists a continuous  $\pi : X \to Y$ mapping $f$ to $g$ on $Y$ with $g$ a
minimal, equicontinuous homeomorphism
on the nontrivial compact metrizable space $Y$.
\end{itemize}
\end{enumerate}
\end{cor}

{\bfseries Proof:} (a): $f$ is chain mixing iff $R_m = X \times X$ by Theorem \ref{theom09}.
On the other hand, $R_m$ is a proper
subset of $X \times X$ iff the quotient space $X/R_m$ is nontrivial in which case $f_m$ on
$X/R_m$ is either a single periodic orbit
(when $X/R_m$ is finite) or an odometer (when $X/R_m$ is infinite) by Theorem \ref{theom11}.
Any odometer in turn projects onto a nontrivial periodic orbit.

The proofs of (b) and (c) use similar applications of Theorem \ref{theom09} and Theorem \ref{theom11}.

$\Box$ \vspace{.5cm}

We immediately obtain the following.

\begin{cor}\label{corm12a} Let $f$ be a closed relation on a compact metric space $(X,d)$.
Assume that $f$ has a fixed point, i.e. \ there
exists $x \in X$ such that $x \in f(x)$. If $f$ chain transitive, strong chain transitive or vague transitive, then
$f$ is chain mixing, strong chain mixing, or vague mixing, respectively. \end{cor}

{\bfseries Proof:} Any factor of $f$ has a fixed point and a nontrivial minimal map does not admit a fixed point.

$\Box$ \vspace{.5cm}

In the case when $f$ is vague transitive we call $R_g$ the \emph{equicontinuous structure relation}.
The induced map $f_g$ on the space
$X/R_g$ is the maximum equicontinuous factor of $f$ on $X$ by Theorem \ref{theom11}(c).

We conclude the section with a question raised by Corollary \ref{corm11a} (b).

\begin{ques}\label{quesm13} Assume $f$ is a continuous map on a compact metric space $(X,d)$
with $\A_d f = X \times X$ and $R_{\ell} = 1_X$.
By Theorem \ref{theom11} and Corollary \ref{corm11a}, $f$ is an equicontinuous homeomorphism
and $d_f \geq d \geq \rho_d^f$ are admissible
metrics on $X$ with $f$ an isometry on $(X,d_f)$ and $(X,\rho_d^f)$.  Are these metrics
necessarily bi-Lipschitz equivalent?
 \end{ques}

The answer is no.

\begin{ex}\label{exm13a} Let $f$ be an irrational rotation on the circle $X$, so that
$f$ is a minimal, equicontinuous homeomorphism. There
exist metrics $d, \rho$ on $X$ with $d \geq \rho$, such that $f$ is an isometry on
$(X,\rho)$ and $(X,d_f)$ but the metrics $d$ and $d_f$  are not
bi-Lipschitz equivalent. \end{ex}

{\bfseries Proof:} We let $X = \R/\Z$ and we can represent it as
$[-\frac{1}{2}, \frac{1}{2}]$ with $-\frac{1}{2}$ identified with $\frac{1}{2}$.
Define $\rho(x_1,x_2) = \inf \{|x_1 - x_2 + n| : n \in \Z \}$. So, of course,
this is a symmetric function of the pair $(x_1 + \Z, x_2 + \Z)$, and
there exists $n_{12} \in \Z$ such that  $\rho(x_1,x_2) = |x_1 - x_2 + n_{12}|$. Hence there exists $n_{23}$ such that
\begin{align}\label{m16}
\begin{split}
 \rho(x_1,x_2) + \rho(x_2,x_3) \ = \ |x_1 - (x_2 - n_{12})| + &|(x_2 - n_{12}) - x_3 + n_{23}| \\
  \geq \ |x_1 - x_3 + n_{23}| \ \geq \ \rho(x_1,x_3)&.
 \end{split}
 \end{align}

 Note that for $x_1, x_2 \in [-\frac{1}{2}, \frac{1}{2}]$, if $ 0 \leq x_1 - x_2 \leq 1$
 then $-1 \leq x_1 - x_2 - 1 \leq 0$ and if
 $ -1 \leq x_1 - x_2 \leq 0$ then $0 \leq x_1 - x_2 + 1 \leq 1$.  It follows in either case that
 \begin{equation}\label{m17}
 x_1, x_2 \in [-\frac{1}{2}, \frac{1}{2}] \quad \Longrightarrow \quad
 \rho(x_1,x_2) \ = \ \min(|x_1 - x_2|, 1 - |x_1 - x_2|).
 \end{equation}

  Thus, $\rho$ is a continuous, translation-invariant metric on $X$, i. e. \ $\rho(x_1, x_2) = \rho(x_1 + c,x_2 + c)$ for
  all $c \in \R/\Z$.  In particular, $f$ is an isometry on $(X,\rho)$.

  For any continuous metric $d$ on $X$, since the iterates of $f$ are dense
  in the group of all translations, it follows that
   \begin{equation}\label{m17a}
   d_f(x_1,x_2) \ = \ \sup \{ d(x_1 + c,x_2 + c) : c \in \R/\Z \}. \hspace{2cm}
   \end{equation}

  For any continuous metric $d$ on $X$, $\A_d f \supset \G f = X \times X$ because $f$ is minimal. If $d \geq \rho$ then
  $\A_{d \times d} (f \times f) \subset \A_{\rho \times \rho} (f \times f)$ and so $R_{\ell} = 1_X$ with respect to $d$.

  On $\R$ we define $q(x) = \sqrt{|x|}$ for $x \in [-\frac{1}{2}, \frac{1}{2}]$
  and extend to a periodic function with period $1$. On
  $X = \R/\Z$ we define
    \begin{equation}\label{m18}
  d(x_1,x_2) \ = \ \max(\rho(x_1,x_2),q(x_1) - q(x_2)|). \hspace{2cm}
  \end{equation}
  So that for $x_1, x_2 \in [-\frac{1}{2}, \frac{1}{2}]/(-\frac{1}{2} = \frac{1}{2})$
  \begin{equation}\label{m18a}
  d(x_1,x_2) \ = \ \max(\min(|x_1 - x_2|, 1 - |x_1 - x_2|),|\sqrt{|x_1|} - \sqrt{|x_2|}|). \hspace{2cm}
  \end{equation}
   Clearly, $d \geq \rho$.

  Since, $ \sqrt{|x_1 - x_2|} \geq  |\sqrt{|x_1|} - \sqrt{|x_2|}| $, i.e.
  $\sqrt{|x_1 - x_2|}$ is a metric on $\R$, it follows that
  if $|x_1 - x_2| < \frac{1}{2}$ the largest value of
  $|\sqrt{|x_1 + t|} - \sqrt{|x_2 + t|}|$ occurs with $t = -x_2$ and so
   \begin{align}\label{m19}
   \begin{split}
  x_1, x_2 \in [-\frac{1}{2}, \frac{1}{2}] \quad &\text{and} \quad |x_1 - x_2| \leq \frac{1}{2} \qquad \Longrightarrow \\
  \rho(x_1,x_2) \ = \ |x_1 - x_2|, \quad &\text{and} \quad d_f(x_1,x_2) \ = \ \max(|x_1 - x_2|,\sqrt{|x_1 - x_2|}).
  \end{split}
  \end{align}

  To see this, observe first that since the square root function is increasing and has a decreasing derivative,
  $(x_1,x_2) \mapsto \sqrt{|x_1 - x_2|}$ is a metric on $\R$. Hence, for $x_1, x_2 \in [-\frac{1}{2}, \frac{1}{2}]$
  \begin{equation}\label{m19a}
  \sqrt{|x_1 - x_2|} \ \geq  \ |\sqrt{|x_1|} - \sqrt{|x_2|}| \ = \ |q(x_1) - q(x_2)|.
  \end{equation}

  Now assume that $0 < x_1 - x_2 \leq \frac{1}{2}$. It follows that
  with $c = -x_2$ $|q(x_1 + c) - q(x_2 + c)| = \sqrt{|x_1 - x_2|}$.

  We now consider translation by an arbitrary  $c \in \R$.

  {\bfseries Case 1} ( $n - \frac{1}{2} \leq x_2 + c < x_1 + c \leq n + \frac{1}{2}$ ):
  In this case,  $x_1 + c - n, x_2 + c - n \in [-\frac{1}{2}, \frac{1}{2}]$
  and so $|q(x_1 + c) - q(x_2 + c)| \leq \sqrt{|x_1 - x_2|}$ by (\ref{m19a}).

  {\bfseries Case 2}  ( $x_2 + c \leq n + \frac{1}{2} < x_1 + c $ ): Since
  $ x_1 - x_2 < \frac{1}{2}$, $x_1 + c - n -1, x_2 + c - n \in [-\frac{1}{2}, \frac{1}{2}]$
   and so $-x_1 - c + n + 1 \in [-\frac{1}{2}, \frac{1}{2}]$. Hence,
      \begin{align}\label{m19b}
      \begin{split}
    |q(x_1 + c) - q(x_2+ c)| \ = \ &|\sqrt{|-x_1 - c + n + 1|} - \sqrt{|x_2 + c - n|}| \\
     \leq \    &\sqrt{|1 - x_1 - x_2 - 2c + 2n|},
    \end{split}
    \end{align}
    by (\ref{m19a}). The assumed inequalities imply that $\pm(1 - x_1 - x_2 - 2c + 2n) \leq x_1 - x_2$ and so
    $|q(x_1 + c) - q(x_2 + c)| \leq \sqrt{|x_1 - x_2|}$.

   {\bfseries Case 3}  ( $x_2 + c < n - \frac{1}{2} \leq   x_1 + c$ ): Since
   $ x_1 - x_2 < \frac{1}{2}$, $x_1 + c - n, x_2 + c - n + 1 \in [-\frac{1}{2}, \frac{1}{2}]$
   and so $-x_2 - c + n - 1 \in [-\frac{1}{2}, \frac{1}{2}]$. Hence,
        \begin{align}\label{m19c}
         \begin{split}
    |q(x_1 + c) - q(x_2+ c)| \ = \ &|\sqrt{|x_1 + c - n|} - \sqrt{|-x_2 - c + n - 1|}| \\
     \leq \    &\sqrt{|-1 - x_1 - x_2 - 2c + 2n|},
    \end{split}
    \end{align}
     by (\ref{m19a}). Again the assumptions imply that $\pm(-1 - x_1 - x_2 - 2c + 2n) \leq x_1 - x_2$ and so
    $|q(x_1 + c) - q(x_2 + c)| \leq \sqrt{|x_1 - x_2|}$.

  This completes the proof of (\ref{m19}).

  Since $\sqrt{|x_1 - x_2|}/|x_1 - x_2| \to \infty$ as $|x_1 - x_2| \to 0$, it
  follows that $d_f$ is not bi-Lipschitz equivalent to $\rho$ on
  any open subset of $X$. On the other hand, $d$ and $\rho$ are bi-Lipschitz
  equivalent on the complement of any closed neighborhood of
  $0 \in X =  [-\frac{1}{2}, \frac{1}{2}]/(-\frac{1}{2} = \frac{1}{2})$,
  because on such a set, the square-root function has a bounded derivative.
  It follows that $d_f$ is not bi-Lipschitz equivalent to $d$ on $X$.

$\Box$ \vspace{.5cm}

We have the following partial result.

\begin{theo}\label{theom14} Assume that $f$ is an homeomorphism on
compact metric space $(X,d)$ with $\A_d f = X \times X$ and $R_{\ell} = 1_X$.
The following are equivalent.
\begin{itemize}
\item[(i)] There exists a continuous metric $d_1$ on $X$ which is bi-Lipschitz
equivalent to $d$ such that $f$ is an isometry on $(X,d_1)$.
\item[(ii)] The metric $d$ is bi-Lipschitz equivalent to $d_f$.
\item[(iii)] The metric $d$ is bi-Lipschitz equivalent to $\rho_d^f$.
\item[(iv)] The metric $d_f$ is bi-Lipschitz equivalent to $\rho_d^f$.
\end{itemize}
\end{theo}

{\bfseries Proof:}  Clearly, (ii), (iii) $\Rightarrow$ (i). Since
$d_f \geq d \geq \rho_d^f$, (iv) $\Rightarrow$ (ii) and (iii).

(i) $\Rightarrow$ (ii), (iii): If for some $K \geq 1$, $(1/K) d \leq d_1 \leq K d$ then,
clearly, $(1/K) d_f \leq (d_1)_f \leq K d_f$ and
$(1/K) (d \times d) \leq (d_1 \times d_1)\leq K (d \times d)$  and so
$(1/K) \ell_{(d \times d)}^{f \times f} \leq \ell_{(d \times d)}^{f \times f}\leq K \ell_{(d \times d)}^{f \times f}$.
Hence,
$(1/K) \rho_d^f \leq \rho_{d_1}^f  \leq K \rho_d^f $.  That is, $d_f$ is bi-Lipschitz equivalent to $(d_1)_f$ and
$\rho_d^f$ is bi-Lipschitz equivalent to $\rho_{d_1}^f$ because  $d$ is bi-Lipschitz equivalent to $d_1$.

Since $f$ is an isometry on $(X,d_1)$, $d_1 = (d_1)_f$ and so $d_1$ is
bi-Lipschitz equivalent to $d_f$. Hence, (i) implies (ii).

In general, $d_1 \geq \rho_{d_1}^f$. On the other hand, since
$f \times f$ is an isometry on $(X \times X, d_1 \times d_1)$,
 Lemma \ref{lemm13} implies that for $x,y \in X$ there exists a point
 $(x_1,x_1) \in \RR (f \times f)(z,z) \subset 1_X$ such that
 $\rho_{d_1}^{f}(x,y) = \ell_{d_1 \times d_1}^{f \times f}((z,z),(x,y)) = \max(d_1(x_1,x), d_1(x_1,y))$.
 Since $d_1(x,y) \leq d_1(x_1,x) + d_1(x_1,y)$,
 it follows that $\rho_{d_1}^{f}(x,y) \geq (1/2) d_1(x,y)$. That is, $d_1$ is
 bi-Lipschitz equivalent to $\rho_{d_1}^{f}$. Hence, (i) implies (iii).

$\Box$ \vspace{1cm}

\section{Transitivity and Weak Mixing for Maps}
\vspace{.5cm}

For a closed relation $f$ on a compact, metrizable space $X$ we defined
$\NN f = \overline{\bigcup_{n \in \N} f^n}$. Since the
twist map $(x,y) \mapsto (y,x)$ is a homeomorphism on $X \times X$ we
see that $\NN (f^{-1}) = (\NN f)^{-1}$ and so we may omit the
parentheses.

The relation $\RR f$ is defined by $\RR f(x) = \overline{\bigcup_{n \in \N} f^n(x)}$
for every $x \in X$. Of course, $\RR f \subset \NN f$
but $\RR f$ is not usually closed and so the inclusion is proper. Furthermore,
$\RR (f^{-1}) $ is usually not equal to $(\RR f)^{-1}$ and so
here the parentheses are required.

When $f$ is a continuous map, Proposition 1.12 of \cite{A93} says that
\begin{align}\label{t01}
\begin{split}
f \cup (f \circ \RR f) \ &= \ \RR f \ = \ f \cup (\RR f \circ f), \\
f \cup (f \circ \NN f) \ &= \ \NN f \ \subset \ f \cup (\NN f \circ f).
\end{split}
\end{align}

\begin{prop}\label{propt01} Let $f$ be a continuous map on a compact
metrizable space $X$.  If $1_X \subset \NN f$, i.e. \ $\NN f$ is reflexive,
then  $\NN f^{-1} = \NN f$, i.e.  $\NN f$ is symmetric. Furthermore, the map $f$ is surjective.\end{prop}

{\bfseries Proof:}  Let $x \in X$. By (\ref{t01}) and induction
$ \NN f \subset f \cup f^2 \cup \dots f^n \cup \NN f \circ f^n$.
Since $(x,x) \in \NN f$, either $f^i(x) = x$ for some $i$ or
$x \in \NN f(f^n(x))$. If $f(x) = x$ then $f^2(x) = x$. If $f^i(x) = x$ with $i \geq 2$ then
then $ f^{(i-1)n}(f^n(x)) = (f^i)^n(x) = x$.  Thus,
in either case, $(f^n(x),x) \in \NN f$. Hence, $\bigcup_{n \in \N} (f^{-1})^n \subset \NN f$.
Since $\NN f$ is closed, $\NN f^{-1} \subset \NN f$.  So $\NN f = (\NN f^{-1})^{-1} \subset \NN f^{-1}$.

From (\ref{t01}) again it follows that either $f(x) = x$ or
there exists $y \in \NN f(x)$ such that $f(y) = x$.  In either case, $x \in f(X)$.
Hence, $f$ is surjective.

$\Box$ \vspace{.5cm}

We call a map $f$ on a compact metrizable space $X$ \emph{topologically transitive}
if $\NN f = X \times X$. There are various, slightly
different notions of topological transitivity in the literature; see
\cite{AC12} where they are sorted out. We follow \cite{A93} where it is shown,
in Theorem 4.12, that $f$ is topologically transitive iff
$Trans_f = \{ x \in X : \RR f(x) = X \}$ is a dense $G_{\d}$ subset of $X$. The map is
minimal iff $Trans_f = X$. The map is called \emph{weak mixing} when
$f \times f$ is topologically transitive on $X \times X$.

\begin{prop}\label{propt02} Let $f$ be a continuous map on a
compact metrizable space $X$. If $\NN f = X \times X$ then for all
$n \in \N$, $\NN f^{(n)}$ is a reflexive, symmetric relation on $X^{(n)}$. \end{prop}

{\bfseries Proof:} By Proposition \ref{propt01} it suffices to
show that $\NN f^{(n)}$ is  reflexive. Let $d$ be a continuous metric on $X$.

Let $\ep > 0$ and $(x_1, \dots, x_n) \in X^{(n)}$. Choose $x \in Trans_f$.
There exist $k_1, \dots, k_n \in \N$ such that $d(f^{k_i}(x),x_i) < \ep/2$
for $i = 1, \dots, n$.
Choose $\d > 0$ so that $o(f^{k_i}, \d) < \ep/2$ for $i = 1, \dots, n$.
There exists $m \in \N$ such that $d(f^m(x),x) < \d$.
Thus, with $(y_1,\dots,y_n) = (f^{k_1}(x),\dots, f^{k_n}(x))$ we have
$d^{(n)}( (f^{(n)})^m(y_1,\dots,y_n),(y_1,\dots,y_n)) < \ep/2$.
Hence,
\begin{align}\label{t02}
\begin{split}
&d^{(n)}((x_1,\dots,x_n),(y_1, \dots, y_n)) < \ep, \\
d^{(n)}&((f^{(n)})^m(y_1,\dots,y_n),(x_1,\dots,x_n)) < \ep.
\end{split}
\end{align}
As $\ep > 0$ was arbitrary, $((x_1,\dots,x_n),(x_1,\dots,x_n)) \in \NN f^{(n)}$.

$\Box$ \vspace{.5cm}

For a continuous map $f$ on a compact space $X$, we call a pair $(x,y) \in X \times X$ \emph
{proximal} when $\RR (f \times f)(x,y) \cap 1_X \not= \emptyset$
and \emph{regionally proximal} when $\NN (f \times f)(x,y) \cap 1_X \not= \emptyset$, thus the \emph{proximality}
and \emph{regional proximality} relations are given by
\begin{equation}\label{t03}
Prox \ = \  (\RR (f \times f))^{-1}(1_X), \qquad Q\ = \  \NN (f \times f)^{-1}(1_X).
\end{equation}
The relations $Prox$ and $Q$ are symmetric and reflexive, but not
usually transitive.  $Q$ is closed, but $Prox$ need not be.

When $\NN f = X \times X$ then Proposition \ref{propt02} implies that
\begin{equation}\label{t04}
Q \ = \ \NN (f \times f)^{-1}(1_X) \ = \ \NN (f \times f)(1_X).
\end{equation}

From (\ref{t01}) it follows that $(f \times f)\circ \NN (f \times f) \subset \NN (f \times f)$.
Hence, $(f \times f)(Q) \subset Q$.
Equivalently, $Q \subset (f \times f)^{-1}(Q)$.

When $\NN f = X \times X$ we define the following relation
\begin{align}\label{t05}
\begin{split}
R_n \ = \ \bigcap_{x \in X} \ \NN (f \times f)^{-1}&(x,x) \ = \ \bigcap_{x \in X} \ \NN (f \times f)(x,x) \ = \\
\{ (x,y) : \NN (f \times f)(x,y) \supset 1_X \} \ &= \ \{ (x,y) : \NN (f \times f)^{-1}(x,y) \supset 1_X \}.
\end{split}
\end{align}

The relation $R_n$ is clearly closed and symmetric.  Since $f$ is
topologically transitive, $1_X \subset \NN (f \times f)(x,x)$ for all
$x \in X$ and so $R_n$ is reflexive as well.

Because $\NN (f \times f)(x,y)$ is closed and $(f \times f)$ $^+$invariant, it follows that
\begin{equation}\label{t08}
R_n \ = \ \NN (f \times f)^{-1}(x^*,x^*) \ = \ \NN (f \times f)(x^*,x^*) \quad \text{for} \ x^* \in Trans_f.
\end{equation}

Recall that  a point $x$ is an \emph{equicontinuity point} when for every
 $\ep > 0$  there exists $\d > 0$
such that $d(x,y) < \d$ implies
$(f^n(x),f^n(y)) \leq \ep$ for $n = 0,1,\dots$.

\begin{lem}\label{lemt03} For a point $x$ the following are equivalent.
\begin{itemize}
\item[(i)] The point $x$ is an equicontinuity point.
\item[(ii)] $ \NN (f \times f)(x,x) = \RR (f \times f)(x,x)$.
\item[(iii)] $ \NN (f \times f)(x,x) \subset 1_X$.
\item[(iv)] If $\{ x_k \}$ and $\{ n_k \}$ are sequences in $X$ and $\N$ with $\{ x_k \} \to x$,
$\{ f^{n_k}(x_k) \} \to y_1$, $\{ f^{n_k}(x) \} \to y_2$, then $y_1 = y_2$.
\end{itemize}
\end{lem}

{\bfseries Proof:} (i) $\Rightarrow$ (ii): Since $(x,x)$ is an
equicontinuity point for $f \times f$, it suffices to show
that $\NN f(x) = \RR f(x)$ when $x$ is an equicontinuity point for $f$.

Suppose $y \in \NN f(x)$, $\ep > 0$, and $\ep > \d > 0$ is chosen as above for the equicontinuity point $x$. There exist
$x_1 \in X$ and $n \in \N$ such that $d(x,x_1) < \d$ and
$d(f^n(x_1),y) < \d$. Then $d(f^n(x),y) \leq d(f^n(x),f^n(x_1)) + d(f^n(x_1),y) < 2 \ep$.

(ii) $\Rightarrow$ (iii), and (iii) $\Rightarrow$ (iv) are obvious.

(iv) $\Rightarrow$ (i): Suppose $x$ is not an equicontinuity point.
There exists $\ep > 0$ so that for every $k \in \N$ there
exists $x_k \in X$ and $n_k \in \N$ with $d(x,x_k) < 1/k$ and
$d(f^{n_k}(x),f^{n_k}(x_k)) \geq \ep$. By going to
a subsequence we can assume $\{ f^{n_k}(x_k) \} \to y_1$ and
$\{ f^{n_k}(x) \} \to y_2$ and we have $d(y_1,y_2) \geq \ep$.

$\Box$ \vspace{.5cm}

If a topologically transitive map admits an equicontinuity point,
then it is a homeomorphism and the set of
equicontinuity points coincides with the residual set of transitive points.
See, e.g. \ \cite{AAuB96} Theorems 2.4, 3.6 and Lemma 3.3.
Such a map is called \emph{almost equicontinuous}. Thus, if a
topologically transitive map is equicontinuous then it is a minimal homeomorphism.
If a minimal map admits an equicontinuity point then it is an equicontinuous homeomorphism.

\begin{theo}\label{theot04} Let $f$ be a topologically transitive map on a compact metrizable space $X$.
\begin{enumerate}
\item[(a)] The map $f$ is equicontinuous iff $Q = 1_X$.

\item[(b)] The map $f$ is almost equicontinuous iff $R_n = 1_X$.
\end{enumerate}
\end{theo}

{\bfseries Proof:} (a): From Lemma \ref{lemt03} it is clear that
every point is an equicontinuity point iff $\NN (f \times f)(1_X) \subset 1_X$.

(b): Similarly, a transitive point $x^*$ is an equicontinuity point iff $\NN (f \times f)(x^*,x^*) \subset 1_X$.

Notice that since $f$ is topologically transitive $\NN (f \times f)(x,x) \supset 1_X$ for any $x \in X$.

$\Box$ \vspace{.5cm}

If $h : X \to Y$ is a continuous map and $x \in X$, the map $h$ is \emph{open at $x$} if $h(U)$ is a
neighborhood of $h(x)$ whenever $U$ is a neighborhood of $x$.

\begin{lem}\label{lemt04a} Let $h : X \to Y$ be a continuous
surjection between compact metrizable spaces.
\begin{enumerate}
\item[(a)] If $h$ is open at $x \in X$ and $\{y_k \}$ is a sequence in $Y$ converging to $h(x)$, then
there exists a sequence $\{x_k \}$ in $X$ converging to $x$ such that $h(x_k) = y_k$ for
all $k$.

\item[(b)] The following conditions are equivalent. When they hold we call $h$ an \emph{almost open map}.
\begin{itemize}
\item[(i)] If $A^{\circ} \not= \emptyset$ for $A \subset X$ then $h(A)^{\circ} \not= \emptyset$.
\item[(ii)] If $U$ is open in $X$ then $U \cap h^{-1}(h(U)^{\circ})$ is an open subset dense in $U$.
\item[(iii)] The set $\{ x \in X : h$ is open at $x \}$ is a dense subset of $X$.
\item[(iv)] The set $\{ x \in X : h$ is open at $x \}$ is a dense, $G_{\d}$ subset of $X$.
\item[(v)] If $D$ is a dense, open subset of $Y$ then $h^{-1}(D)$ is a dense open subset of $X$.
\item[(vi)] If $D$ is a dense subset of $Y$ then $h^{-1}(D)$ is a dense subset of $X$.
\end{itemize}
\end{enumerate}
\end{lem}

{\bfseries Proof:} (a): Given a metric $d $  on $X$, there is an increasing sequence $N_0 < N_1 < \dots$
such that $y_k \in h(V_{1/k}(x))$ for $k \geq N_k$. So we can choose $x_k \in X$ such that
$h(x_k) = y_k$ and $x_k \in V_{1/k}(x)$ for $k \geq N_k$.

(b): It is clear that (vi) $\Rightarrow$ (v), and (iv) $\Rightarrow$ (iii) $\Rightarrow$ (ii) $\Rightarrow$ (i).

(i) $\Rightarrow$ (ii): If $V$ is an open subset of $U$ then (i) implies that $V \cap h^{-1}(h(V)^{\circ})$
is nonempty and is contained in $V \cap (U \cap h^{-1}(h(U)^{\circ}))$.

(ii)  $\Rightarrow$ (iv): Let $\A_k$ be a finite cover of $X$  by open sets of diameter less than $1/k$.
Let $U_{k} = \bigcup_{A \in \A_k} \ A \cap h^{-1}(h(A)^{\circ})$.
Each $U_{k}$ is open and dense and so by the Baire Category Theorem, $D = \bigcap_{k} U_{k}$ is
a dense, $G_{\d}$ subset.  It is easy to check that $D$ is the set of points at which $h$ is open.

(i) $\Rightarrow$ (vi): If $D$ is dense in $Y$ and $U$ is open in $X$ then by (i) $h(U)$ meets $D$ and
so $U$ meets $h^{-1}(D)$.

(v) $\Rightarrow$ (i): If $A^{\circ} \not= \emptyset $ but $h(A)^{\circ} = \emptyset$, then we can
choose $C$ a closed subset of  $A^{\circ}$ with $C^{\circ} \not= \emptyset $. Hence, $h(C)$ is
compact with $h(C)^{\circ} = \emptyset $.  Hence, $D = Y \setminus h(C)$ is open and dense, but
$h^{-1}(D) \cap C = \emptyset$ and so $h^{-1}(D)$ is not dense.

$\Box$ \vspace{.5cm}

Let $f$ on $X$ and $g$ on $Y$ be continuous surjections of compact metrizable spaces. If $h : X \to Y$ is a continuous
surjection which maps $f$ to $g$ on $Y$, then we call $h$ a \emph{minimal morphism} if $X$ is the only closed, $f$ invariant subset
which is mapped by $h$ onto $Y$.

\begin{lem}\label{lemt04b} Let $h : X \to Y$ be a continuous surjection mapping $f$ on $X$ to $g$ on $Y$.
\begin{enumerate}
\item[(a)] If $h$ is a minimal morphism and $A$ is a closed, $f$ $^+$invariant subset of $X$ with $h(A) = Y$, then $A = X$.

\item[(b)] If $f$ is topologically transitive, and so $g$ is topologically transitive, then the following are equivalent.
\begin{itemize}
\item[(i)] The map $h$ is a minimal morphism.
\item[(ii)] There exists a transitive point $y$ for $g$ such that every $x \in h^{-1}(y)$ is a transitive point for $f$.
\item[(iii)]For every transitive point $y$ for $g$, every $x \in h^{-1}(y)$ is a transitive point for $f$.
\end{itemize}
\end{enumerate}
\end{lem}

{\bfseries Proof:} (a): The intersection $B = \bigcap_{n \in \N} f^n(A)$ is a closed, $^+$invariant subset and if
$x \in B$ then $\{ f^{-1}(x) \cap f^n(A) \}$ is a non-increasing sequence of nonempty compacta.  The intersection is
in $f^{-1}(x) \cap B$ and so $B$ is $f$ invariant. For each $n \in \N$, $h(f^n(A)) = g^n(h(A)) = g^n(Y) = Y$ since $g$
is surjective. So for $y \in Y$, $\{ h^{-1}(y) \cap f^n(A) \}$ is a non-increasing sequence of nonempty compacta with intersection is
$h^{-1}(y) \cap B$. That is, $h(B) = Y$. Since $h$ is minimal, $X = B \subset A$.

(b): A factor of a topologically transitive map is topologically transitive and if $x$ is a transitive point for $X$ then
$h(x)$ is a transitive point for $g$ since $h(\RR f(x)) = \RR g(h(x))$.

(i) $\Rightarrow$ (iii): Let $h(x) = y$ be a transitive point for $g$.  Then $h(\RR f(x)) = \RR g(y) = Y$. Since $\RR f(x)$ is
$f$ $^+$invariant, it follows from (a) that $\RR f(x) = X$, i.e. $x$ is a transitive point for $f$.

(iii) $\Rightarrow$ (ii): Obvious since $Y$ has transitive points.

(ii) $\Rightarrow$ (i): If  $h(A) = Y$ then $y \in h(A)$ and so there exists $x \in h^{-1}(y) \cap A$.
By assumption on $y$, $x$ is a transitive point for $f$. If $A$ is closed and invariant, then
$X = \RR f(x) \subset A$.

$\Box$ \vspace{.5cm}

\begin{prop}\label{propt04c} Let $f$ be a topologically transitive
map on a compact metrizable space $X$ and let $h : X \to Y$ be a continuous surjection mapping $f$ on $X$ to $g$ on $Y$
with $Y$ metrizable.
\begin{enumerate}
\item[(a)] If $h$ is almost open or a minimal morphism, then there exists a transitive point $x$ for $f$ such that
$h$ is open at $x$.

\item[(b)] If $g$ is almost open, e.g. \ if $g$ is a homeomorphism, and there exists a transitive point $x$ for $f$ such that
$h$ is open at $x$, then $h$ is almost open. In particular, if $h$ is a minimal morphism and $g$ is almost open, then $h$ is
almost open.
\end{enumerate}
\end{prop}

{\bfseries Proof:} (a): If $h$ is almost open, then the set of points and which $h$ is  a dense $G_{\d}$, as is the set of
transitive points for $f$. By the Baire Category Theorem, the intersection is nonempty and these are transitive points of $f$ at
which $h$ is open.

For the minimal morphism case, we quote some results from \cite{A93}. The relation $h^{-1} : Y \to X$ is closed and so is
upper semicontinuous. Let $D \subset Y$ be the set of points $y$ at which $h^{-1}$ is lower  semicontinuous.  These are exactly
the set of points $y$ such that $h$ is open at every point of $h^{-1}(y)$.
By Theorem 7.19 of \cite{A93} the set $D$ is a dense $G_{\d}$ and so contains a transitive point $y$. So $h$ is open at every
point $x$ of $h^{-1}(y)$ and when $h$ is a minimal morphism, these are transitive points for $f$.

(b): For every $n \in \N$, $g^n \circ h = h \circ f^n$ and so $g^n \circ h \circ f^{-n} = h \circ f^n \circ f^{-n} = h$, since $f$ is surjective.
Now let $U$ be open in $X$.  There exists $n \in \N$ such that $f^n(x) \in U$ and so $f^{-n}(U)$ is an open set containing $x$. Since
$h$ is open at $x$ and $g$ is almost open, $h(U) = g^n(h(f^{-n}(U)))$ has a nonempty interior.

$\Box$ \vspace{.5cm}

For the next result we repeat and adapt the lovely proof of Auslander's second folk theorem from \cite{Au16}.

\begin{theo}\label{theot05} Let $f$ be a topologically transitive
map on a compact metrizable space $X$ and $h : X \to Y$ be a continuous
surjection with $Y$ metrizable. Assume that $h$ maps $f$ to $g$ a
continuous map on $Y$.
\begin{enumerate}
\item[(a)] If $Q \subset (h \times h)^{-1}(1_Y)$ then
$g$ is equicontinuous on $Y$.

\item[(b)] If there exists a transitive point $x$ for $f$ such that
$h$ is open at $x$, and $R_n \subset (h \times h)^{-1}(1_Y)$ then
$g$ is almost equicontinuous on $Y$ and $h$ is almost open.
\end{enumerate}
\end{theo}

{\bfseries Proof:} (a): Let $\{ y_k \}$ and $\{ n_k \}$ be sequences in $Y$ and $\N$ with $\{ y_k \} \to y$,
$\{ g^{n_k}(y_k) \} \to z_1$, $\{ g^{n_k}(y) \} \to z_2$. There exist $x_k \in X$ with $h(x_k) = y_k$ and by going to
a subsequence we can assume that  $\{ x_k \} \to x$, $\{ f^{n_k}(x_k) \} \to w_1$ and $\{ f^{n_k}(x) \} \to w_2$.
Hence, $h(x) = y, h(w_1) = z_1$ and $h(w_2) = z_2$. Since
$\{ (x_k,x) \} \to (x,x)$ and $(f^{n_k}(x_k),f^{n_k}(x)) \} \to (w_1,w_2)$,
it follows that $(w_1,w_2) \in Q \subset (h \times h)^{-1}(1_Y)$. Hence, $y_1 = h(w_1) = h(w_2) = y_2$.  Thus, $y$ is an
equicontinuity point for $g$. Since $y$ was arbitrary, $g$ is equicontinuous.

(b): Let $x$ be a transitive point  for $f$ at which $h$ is open.  Let $y = h(x)$.
In the above proof we can choose $x_k$ so that $\{ x_k \} \to x$, by Lemma \ref{lemt04a} (a). In the above
proof we obtain that  $(w_1,w_2) \in R_n \subset (h \times h)^{-1}(1_Y)$ since $x$ is a transitive point.
As before, $y_1 = h(w_1) = h(w_2) = y_2$.  Thus, $y$ is an
equicontinuity point for $g$. Since $g$ is transitive, it follows that $g$ is almost equicontinuous. Since a topologically transitive
almost equicontinuous map is a homeomorphism, it follows from Proposition \ref{propt04c} (b) that $h$ is almost open.

$\Box$ \vspace{.5cm}

In particular, this shows that a factor of an equicontinuous
minimal homeomorphism is an equicontinuous minimal homeomorphism.  Because an
almost equicontinuous topologically transitive map can admit
factors which are not almost equicontinuous, the analogous result is not true for
$R_n$ and general $h$.  A factor of an almost equicontinuous, transitive map by an
almost open mapping is almost equicontinuous. This is Lemma 1.6 of \cite{GW}.

\begin{cor}\label{cort06} If $f$ is a topologically transitive
map on a compact metrizable space $X$, then $\G Q = R_g$.  That is,
$R_g$ is the smallest closed transitive relation which contains $Q$. \end{cor}

{\bfseries Proof:} Clearly, the closed equivalence relation $R_g$ contains $Q$ and so contains $\G Q$.
On the other hand, if $E = \G Q$ then because
$Q$ is reflexive and symmetric, $E$ is an equivalence relation. Since $Q \subset (f \times f)^{-1}(Q)$ it follows that
$Q $ is contained in the closed equivalence relation $E \cap (f \times f)^{-1}(E)$ and so the latter equals $E$.
Thus, $(f \times f)(E) \subset E$. Hence, if $\bar X = X/E$ then $f$ induces a continuous map $\bar f$ on $\bar X$. Since
$Q \subset E$,  Auslander's second folk theorem \ref{theot05}(a) implies that $\bar f$ is equicontinuous.

Since $f$ is topologically transitive it is vague transitive and so
Theorem \ref{theom11} (c) implies that $R_g \subset E$.

$\Box$ \vspace{.5cm}

For the following analogue of  Theorem \ref{theom09} we apply a number of well-known results.

\begin{theo}\label{theot07} Let $f$ be a continuous map on a compact
metrizable space. If $f$ is topologically transitive then
the following are equivalent.
\begin{itemize}
\item[ (i)]  The map $f$ is weak mixing.
\item[(ii)] The map $f^{(n)}$ on $X^{(n)}$ is weak mixing for every $n \in \N$.
\item[(iii)] $R_n = X \times X$.
\end{itemize}
\end{theo}

{\bfseries Proof:} (ii) $\Rightarrow$ (i) $\Rightarrow$ (iii) are obvious.

(iii) $\Rightarrow$ (i) and (i) $\Rightarrow$ (ii) are applications of the
sharpening by Karl Peterson of the beautiful Furstenberg
Intersection Lemma, see, e.g.  \ \cite{A97} Lemma 4.1 and Proposition 4.10.
In particular, $f$ is weak mixing if for every $x,y \in X$
$(x,y) \in \NN (f \times f)(x,x)$.

$\Box$ \vspace{.5cm}

Thus,  a topologically transitive map  $f$ is weak mixing
iff $R_n = X \times X$.  By Theorem \ref{theom11} (c)
the map $f$ on $X$ has a non-trivial equicontinuous factor
iff $R_g$ is a proper subset of $X \times X$. The analogue of the
dichotomy result, Corollary \ref{corm12}, fails because of the gaps between $R_n$ and $Q$ and $R_g$.

\begin{theo}\label{theot08} Let $f$ be a topologically transitive, continuous map on a compact metrizable space $X$.
\begin{enumerate}
\item[(a)] The following are equivalent.
\begin{itemize}
\item[(i)] The relation $Q$ is transitive and so is a closed equivalence relation.
\item[(ii)] $Q = R_g$.
\item[(iii)] $\NN (f \times f)(Q) = Q$.
\end{itemize}

\item[(b)] If $\NN (f \times f)(R_n) = R_n$ then $R_n = Q = \NN (f \times f)(Q)$.
\end{enumerate}
\end{theo}

{\bfseries Proof:} (a): (i) $\Rightarrow$ (ii): If $Q$ is transitive,
then $Q = \G Q$ which equals $R_g$ by Corollary \ref{cort06}.

(ii) $\Rightarrow$ (iii): Since $\G (f \times f)$ is transitive,
\begin{align}\label{t09}
\begin{split}
\NN (f \times f)(Q) \ \subset  \ &\G (f \times f)(R_g) \ = \\
\G (f \times f)\G (f \times f)(1_X) \ \subset \ &\G (f \times f)(1_X) \ = \ R_g \ = \ Q.
\end{split}
\end{align}

(iii)  $\Rightarrow$ (i): This result, and its proof, come from \cite{AuG97}
Lemma 5. Assume that $(x,y), (y,z) \in Q$. Since $(y,z) \in Q$,
there exist $y_k,z_k, w \in X$ and $n_k \in \N$ such that
$\{ (y_k,z_k) \} \to (y,z)$ and $ \{ (f^{n_k}(y_k),f^{n_k}(y_k)) \} \to ( w,w)$.
We can assume that $\{ f^n_k(x) \} \to w_1 \in X$. Now
$\{ (x,y_k) \} \to (x,y) \in Q$ and $\{ (f^{n_k}(x),f^{n_k}(y_k)) \} \to (w_1,w)$
and so $(w_1,w) \in \NN (f \times f)(Q) = Q$. Also
$\{ (x,z_k) \} \to (x,z) $ and $\{ (f^{n_k}(x),f^{n_k}(z_k)) \} \to (w_1,w) \in Q$.
Hence $(x,z) \in \NN (f \times f)^{-1}(Q)$ which equals $Q$ since $\NN (f \times f)$ is symmetric.

(b): Since $1_X \subset R_n \subset Q = \NN (f \times f)(1_X)$
we have $R_n \subset Q \subset \NN (f \times f)(R_n)$. So if
$R_n = \NN (f \times f)(R_n)$, then $R_n = Q $ and $Q = R_n =  \NN (f \times f)(R_n) =  \NN (f \times f)(Q)$.

$\Box$ \vspace{.5cm}

If $f$ is a minimal map then every point is transitive and so
$Q = R_n$. It is a much deeper result that for a minimal homeomorphism
$f$, the regional proximality relation $Q$ is transitive and so $Q = R_g$.
This follows from Theorem 8 of Chapter 9 in \cite{Au88} and
is proved directly in \cite{AuG97}. In particular, this yields:

\begin{theo}\label{theot09} Assume that $f$ is a minimal homeomorphism. Exactly one of the following is true:
\begin{itemize}
\item[(i)] The relation $f$ is weak mixing.
\item[(ii)] There exists a continuous  $\pi : X \to Y$ mapping $f$ to
$g$ on $Y$ with $g$ a minimal, equicontinuous homeomorphism
on the nontrivial compact metrizable space $Y$.
\end{itemize}
\end{theo}

$\Box$ \vspace{.5cm}

\begin{cor}\label{cort09a} If $f$ is a minimal homeomorphism, then $f$ is weak mixing if and only if it is vague mixing. \end{cor}

{\bfseries Proof:} Because $f$ is minimal, $R_n = Q = R_g$. Weak mixing is equivalent to $R_n = X \times X$ and vague mixing is
equivalent to $R_g = X \times X$.

$\Box$ \vspace{.5cm}  

Without minimality we can run into the following, compare Corollary \ref{corm12a}.

\begin{ex}\label{ext10} There exists a topologically transitive
homeomorphism $f$ on a compact metrizable space $X$, which admits a fixed point as its unique minimal subset, but is
not weak mixing and so $R_n$ is a proper subset of $X \times X$. On the other hand, $Prox = Q = X \times X$ and 
$R_n \circ R_n = X \times X$ and so  $X \times X$ is the
smallest transitive relation which contains $R_n$.\end{ex}

{\bfseries Proof:} We sketch the construction. Let $f_0$ be a
non-trivial weak mixing homeomorphism on $X_0$ which has a fixed point $e_0$ as unique minimal subset. The
\emph{stopped torus} of \cite{A93} Chapter 9 is one such example.

Let $f_1$ be an irrational rotation on a circle $X_1$ so that
$f_1$ is minimal and equicontinuous. Let $f_2 = f_0 \times f_1$ on
$X_0 \times X_1$.  Because $f_0$ is weak mixing and $f_1$ is
minimal, it follows that $f_2$ is topologically transitive.  Since it has
$f_1$ as a factor, it is not weak mixing.  Now let $\pi : X_0 \times X_1 \to X$ be the
quotient map obtained by identifying the invariant subset $\{ e_0 \} \times X_1$
to a point $e$ and let $f$ be the induced homeomorphism on $X$. Since $\{ e_0 \} \times X_1$ is
the unique minimal subset of $X_0 \times X_1$, the fixed point $e$ is the unique minimal subset of $X$.
As $f$ is a factor of $f_2$ it is topologically transitive. On the other hand, it cannot be weak mixing
because $f_2$ is an almost one-to-one lift and such lifts preserve weak mixing.

The point $(e,e)$ is the unique minimal subset of $X \times X$ and so  $X \times X = (\RR (f \times f))^{-1}(e,e)$. Hence,
$ X \times X = Prox \subset Q \subset X \times X$.

In addition, since $f_0 \times f_0$ is weak mixing on $X_0 \times X_0$, the homeomorphism
$f_0 \times f_0 \times f_1$ is topologically transitive. So if $(x,t), (y,s) \in X_0 \times X_1$
there exists a sequence $\{(u_k,v_k,t_k) \in X_0 \times X_0 \times X_1 \}$ which converges
to $(x,x,t)$ and $n_k \to \infty$ such that $\{(f_0^{n_k}(u_k),f_0^{n_k}(v_k),f_1^{n_k}(t_k)) \}$
converges to $(y,e_0,s)$. We have that $(\pi(u_k,t_k),\pi(v_k,t_k))$ converges to $(\pi(x,t),\pi(x,t))$
and $f^{n_k}(\pi(u_k,t_k)),f^{n_k}(\pi(v_k,t_k))$ converges to \\ $(\pi(y,s),\pi(e_0,s)) = (\pi(y,s),e)$.
Since $x,y,t,s$ are arbitrary, it follows that for all $y_1,y_2 \in X_0, s_1,s_2 \in X_1$, 
$(\pi(y_1,s_1),e) \in R_n$ and $(\pi(y_2,s_2),e) \in R_n$.  Since $R_n$ is symmetric, we have that \\
$(\pi(y_1,s_1),\pi(y_2,s_2)) \in R_n \circ R_n$ and so $R_n \circ R_n = X \times X$.

$\Box$ \vspace{1cm}

\bibliographystyle{amsplain}

\end{document}